\documentclass[12pt, reqno]{amsart}

\usepackage{amssymb, amsfonts, amsmath, amsthm} % standard AMS packages
\newtheorem{introtheo}{Theorem}

\usepackage{mathrsfs}
\usepackage{eucal}
\usepackage{txfonts}
\usepackage{graphicx}
\usepackage[dvipsnames]{xcolor}
\usepackage{tikz}
\usepackage[all]{xy}
\usepackage{multicol}
\usepackage{xspace}
\usepackage{enumitem}
\usepackage{geometry}
\usepackage{bm} % For bold math symbols if needed
\usepackage{dsfont}
\usepackage[
    colorlinks=true,
    linkcolor=blue,
    citecolor=blue,
    urlcolor=blue,
    filecolor=blue,
    final,
    backref=page,
    hyperindex
]{hyperref}

\newtheoremstyle{reference} %
  {}          % Space above
  {}          % Space below
  {\normalfont} % Body font
  {}          % Indent amount
  {\bfseries} % Head font
  {.}         % Punctuation after head
  { .5em}     % Space after head
  {\thmname{#1}\thmnumber{ #2}\thmnote{ \textnormal{(#3)}}} %

\DeclareMathOperator{\End}{End}

\newcommand{\Aut}{\mathrm{Aut}}
\newcommand{\Der}{\mathrm{Der}}

\newcommand{\id}{\mathrm{id}}
\numberwithin{equation}{section}

\newtheorem{theo}{Theorem}[section]
\newtheorem{defi}[theo]{Definition}
\newtheorem{coro}[theo]{Corollary}
\newtheorem{lemm}[theo]{Lemma}
\newtheorem{prop}[theo]{Proposition}

\newtheorem{rema}[theo]{Remark}

\title{New simple modules over the $N=2$ superconformal algebra from Weyl modules}

\author{Haibo Chen}
\address{School of Science, Jimei University, Xiamen, Fujian 361021, China}
\email{hypo1025@jmu.edu.cn}

\author{Dong Liu}
\address{Department of Mathematics, Huzhou University, Zhejiang Huzhou, 313000, China}
\email{liudong@zjhu.edu.cn}

\author[Yufeng Pei]{Yufeng Pei}
\address{Department of Mathematics, Huzhou University, Zhejiang Huzhou, 313000, China}
\email{pei@zjhu.edu.cn}

\begin{document}
\begin{abstract}
Let $\mathfrak{g}$ be the centerless $N=2$ superconformal algebra in its
Ramond form, and let $\mathcal{D}$ be the Weyl algebra generated by
$t^{\pm1}$ and $D$. For each $\alpha,\beta\in\mathbb{C}$ we construct a functor
$F_{\alpha,\beta}$ from $\mathcal{D}$-modules to $\mathfrak{g}$-modules by
twisting the contact realization of $\mathfrak{g}$ with a $1$-cocycle. It
preserves simplicity if and only if $(2\alpha-\beta)(2\alpha+\beta)\neq0$. We
also give a necessary and sufficient condition for two of these simple modules
to be isomorphic. Varying the $\mathcal{D}$-module gives both non-weight
modules and uniformly bounded weight modules with four-dimensional weight
spaces, related by Nilsson's weighting functor.
\end{abstract}

\maketitle
%\tableofcontents

\section{Introduction}
The $N=2$ superconformal algebra is one of the fundamental extended
superconformal algebras in two-dimensional superconformal field theory and its
topological twist \cite{EY90,Gan97,Ki,LVW,Wit88,Vafa91}. As a Lie
superalgebra it belongs to the class studied in \cite{Kac77,KL,SS}. It arises
in several sectors, of which the Neveu--Schwarz and Ramond sectors are
isomorphic \cite{SS}. We work throughout with the centerless Ramond sector
$\mathfrak{g}$, which admits differential operator realizations coming from
contact geometry and super-Riemann surfaces \cite{Ba,KL,RSV}. These
realizations are the starting point for our construction.
By \cite{LPX0}, the simple uniformly bounded weight modules of the $N=2$
Ramond algebra can have four-dimensional weight spaces. The two-component
construction of \cite{CDLP2} does not realize such modules. Producing them
motivates this paper. We obtain a much larger family: for each
$\alpha,\beta\in\mathbb{C}$ and each simple $\mathcal{D}$-module $P$, a
four-component $\mathfrak{g}$-module $F_{\alpha,\beta}(P)$. As $P$ varies, this
family contains both the motivating uniformly bounded weight modules and a
class of $U(\mathbb{C}L_0)$-free non-weight modules.

By \cite{LPX0}, the simple uniformly bounded weight modules of the $N=2$ Ramond
algebra can have four-dimensional weight spaces, which the two-component
construction of \cite{CDLP2} does not realize. Producing such modules motivates
this paper, and we obtain a much larger family: for each $\alpha,\beta\in\mathbb{C}$
and each simple $\mathcal{D}$-module $P$, a four-component $\mathfrak{g}$-module
$F_{\alpha,\beta}(P)$. As $P$ varies, this family contains both the motivating
uniformly bounded weight modules and a class of $U(\mathbb{C}L_0)$-free non-weight
modules.

For any $\mathcal{D}$-module $P$, we build $F_{\alpha,\beta}(P)$ as follows:
\[
\xymatrix@C=3em@R=2.4em{
  P \ar[r]^-{\mathrm{Morita}}
  & P\otimes\Lambda(\theta_1,\theta_2)
    \ar[r]^-{\pi^{\mathrm c}}
  & \text{contact }\mathfrak{g}\text{-module}
    \ar[d]^-{\mathrm{Ad}_\Phi} \\
  && \text{twisted-contact }\mathfrak{g}\text{-module}
    \ar[r]^-{d_{\alpha,\beta}}
  & F_{\alpha,\beta}(P).
}
\]
The Morita equivalence between $\mathcal{D}$ and the two-odd-variable Weyl
superalgebra $\mathcal{SD}(1,2)$, induced by
$\mathcal{SD}(1,2)\cong\mathcal{D}\otimes\mathrm{End}_{\mathbb{C}}
\Lambda(\theta_1,\theta_2)$, sends $P$ to $P\otimes\Lambda(\theta_1,\theta_2)$.
The standard contact realization $\pi^{\mathrm c}$ of $\mathfrak{g}$ on
$A=\mathbb{C}[t^{\pm1},\theta_1,\theta_2]$ \cite{KL,RSV} makes this a
$\mathfrak{g}$-module. The coordinate automorphism $\Phi$ then rewrites the
action so that the odd generators act triangularly on the four summands $P$,
$P\theta_1$, $P\theta_2$, $P\theta_1\theta_2$. Finally, we deform by a
multiplication-valued $1$-cocycle
$d_{\alpha,\beta}\colon\mathfrak{g}\to A$. Such cocycles modulo coboundaries
form a three-dimensional space, but a shift in $D$ absorbs one direction,
leaving $\alpha,\beta$ as the two essential parameters.

Two features of our approach are worth noting. First, we work with the contact
realization on two odd variables, rather than the topological-twist
realization of the one-variable case \cite{CDLP2}. In fact we use a twisted
form of it, given by the coordinate automorphism $\Phi$ above
(Proposition~\ref{prop:contact-asym-conjugacy}). This twisted realization
triangularizes the action on the four summands $P$, $P\theta_1$, $P\theta_2$,
$P\theta_1\theta_2$, and the topological-twist and contact modules are
non-isomorphic (Proposition~\ref{prop:topological-comparison}). Second, we
improve the twisting technique of L\"u and Zhao by working directly with the
$1$-cocycle, following the bifunctor viewpoint of \cite{PSTZ}. Earlier
constructions insert the parameters by hand through an automorphism. In ours
the parameters are the cohomology classes we classify
(Proposition~\ref{prop:multiplication-cocycles}), which makes their meaning
explicit.

We can now state our main results. The first gives a simplicity criterion for
$F_{\alpha,\beta}$, together with the structure of $F_{\alpha,\beta}(P)$ when it
fails.
\begin{introtheo}\label{ithm:simplicity}
Let $P$ be a nonzero simple $\mathcal{D}$-module.
\begin{enumerate}
\item[\textup{(i)}]
The $\mathfrak{g}$-module $F_{\alpha,\beta}(P)$ is simple if and only if
$(2\alpha-\beta)(2\alpha+\beta)\neq0$
\textup{(Theorem~\ref{thm:irreducibility})}.
\item[\textup{(ii)}]
On each singular hyperplane $2\alpha\pm\beta=0$, the module
$F_{\alpha,\beta}(P)$ is an extension of two two-component modules of the kind
constructed in \cite{CDLP2}
\textup{(Proposition~\ref{prop:singular-exact-sequences})}.
\end{enumerate}
\end{introtheo}
The second classifies the simple modules up to isomorphism.
\begin{introtheo}\label{ithm:iso}
Suppose $(2\alpha_i-\beta_i)(2\alpha_i+\beta_i)\neq0$ for $i=1,2$, and let
$P_1,P_2$ be simple $\mathcal{D}$-modules. Then
$F_{\alpha_1,\beta_1}(P_1)\cong F_{\alpha_2,\beta_2}(P_2)$ as
$\mathfrak{g}$-modules if and only if $\alpha_1=\alpha_2$, $\beta_1=\beta_2$,
and $P_1\cong P_2$ \textup{(Theorem~\ref{thm:isomorphism})}.
\end{introtheo}

In particular, $F_{\alpha,\beta}$ preserves simplicity exactly off the two
singular hyperplanes. On them it recovers the two-component modules of
\cite{CDLP2}; off them it produces new simple $\mathfrak{g}$-modules, classified
by Theorem~\ref{ithm:iso}. The simplicity criterion is proved using a finite
difference operator $\mathcal{S}_{k,d}\in U(\mathfrak{g})$. This operator
annihilates $P\theta_1$, $P\theta_2$, $P\Theta$ and maps $P$ into $P\Theta$.
The construction applies to every simple $\mathcal{D}$-module, hence to all
families in the classification \cite{B,LZ}. We treat two in detail. From the
polynomial $\mathcal{D}$-modules $\Omega(\lambda)$ we construct and characterize
a new class of $\mathfrak{g}$-modules, free of rank four over
$U(\mathbb{C}L_0)$ (Theorem~\ref{thm:characterization-standard-rank-four}). The
rank-one weight $\mathcal{D}$-modules $\mathcal{V}_a$ yield the uniformly
bounded weight modules, with four-dimensional weight spaces. Nilsson's weighting
functor with respect to $L_0$ relates the two: it sends the rank-four free
modules to the uniformly bounded ones
(Theorem~\ref{thm:weighting-reduction}).

The paper is organized as follows. Section~\ref{sec:preliminaries} fixes
conventions, records the spectral flow, and reviews the needed facts on Weyl
superalgebras. Section~\ref{sec:Weyl_reps} compares the topological-twist,
contact, and twisted-contact realizations. Section~\ref{sec:twisting-functors}
computes the relevant cohomology and defines $F_{\alpha,\beta}$.
Section~\ref{sec:irreducibility} proves the simplicity criterion, the singular
exact sequences, and the isomorphism theorem. Section~\ref{sec:applications}
treats the two examples and relates them by Nilsson's weighting functor.

\smallskip
\noindent\textbf{Notation.} Throughout, $\mathbb{Z}$, $\mathbb{Z}^*$,
$\mathbb{N}$, and $\mathbb{Z}_+$ denote the integers, nonzero integers,
nonnegative integers, and positive integers, respectively, and $\mathbb{C}$,
$\mathbb{C}^*$ the complex numbers and the nonzero complex numbers. All vector
spaces, algebras, and tensor products are over $\mathbb{C}$. For a Lie
superalgebra $L$, $U(L)$ denotes its universal enveloping algebra; for an
associative algebra $S$, $S\text{-}\mathrm{Mod}$ denotes the category of left
$S$-modules. Submodules are understood to be graded and homomorphisms to be
even; for purely even algebras this agrees with the usual convention. A module
is \emph{simple} if it is nonzero and has no proper nonzero submodules.
\section{Preliminaries}
\label{sec:preliminaries}
This section fixes notation and recalls the background needed below. We
introduce the centerless \(N=2\) superconformal algebra, record the spectral
flow relating its Ramond and Neveu--Schwarz forms, and review the Weyl
superalgebra together with the classification of its simple modules via Morita
equivalence.

We recall the centerless \(N=2\) superconformal algebra.
For \(\epsilon\in\{0,\tfrac12\}\), let \(\mathfrak g[\epsilon]\) be the Lie
superalgebra with even basis \(\{L_m,H_m\mid m\in\mathbb Z\}\) and odd basis
\(\{G_r^+,G_r^-\mid r\in\mathbb Z+\epsilon\}\), subject to the defining
relations
\begin{align}
    [L_m,L_n]      &= (n-m)L_{m+n},                          \label{eq:LL}\\
    [L_m,H_n]      &= nH_{m+n},                              \label{eq:LH}\\
    [L_m,G_r^\pm]  &= \left(r-\tfrac{m}{2}\right)G_{m+r}^\pm,\label{eq:LG}\\
    [H_m,G_r^\pm]  &= \pm G_{m+r}^\pm,                       \label{eq:HG}\\
    \{G_r^+,G_s^-\}&= -2L_{r+s}+(r-s)H_{r+s},                \label{eq:GG_pm}\\
    \{G_r^\pm,G_s^\pm\}&= [H_m,H_n]=0.                       \label{eq:GG_same}
\end{align}
Here \([\,\cdot\,,\,\cdot\,]\) is the superbracket, and
\(\{\,\cdot\,,\,\cdot\,\}\) denotes the bracket when both arguments are odd. All
remaining brackets are determined by super skew-symmetry.
For \(\epsilon=0\) this is the centerless Ramond algebra
\(\mathfrak g:=\mathfrak g[0]\), and for \(\epsilon=\tfrac12\) the
Neveu--Schwarz algebra \(\mathfrak g_{\mathrm{NS}}:=\mathfrak g[\tfrac12]\).

%\begin{rema}
%\label{rem:virasoro-sign}
%Our convention \([L_m,L_n]=(n-m)L_{m+n}\) has the opposite sign from the one
%commonly used in physics. We adopt it because
%\[
%  [t^mD,t^nD]=(n-m)t^{m+n}D .
%\]
%\end{rema}

The Ramond and Neveu--Schwarz presentations are related by half-integral
spectral flow.

\begin{prop}%[Spectral flow]
\label{prop:spectral-flow}
Let \(\eta\in\tfrac12\mathbb Z\), and identify
\(\mathfrak g[\epsilon+k]=\mathfrak g[\epsilon]\) for \(k\in\mathbb Z\).
The linear map
\[
  \sigma_\eta:\mathfrak g[\epsilon]\longrightarrow
  \mathfrak g[\epsilon-\eta]
\]
defined by
\[
  \sigma_\eta(L_m)=L_m+\eta H_m,\quad
  \sigma_\eta(H_m)=H_m,
\]
and
\[
  \sigma_\eta(G_r^+)=G_{r-\eta}^+,\quad
  \sigma_\eta(G_r^-)=G_{r+\eta}^-,
\]
is an isomorphism of Lie superalgebras. In particular, \(\sigma_{1/2}\)
identifies the Ramond and Neveu--Schwarz presentations.
\end{prop}

\begin{proof}
Since \(2\eta\in\mathbb Z\), we have
\(\mathbb Z+\epsilon+\eta=\mathbb Z+\epsilon-\eta\). Hence the shifted indices
\(r-\eta\) and \(r+\eta\) lie in the target index set, and \(\sigma_\eta\) is
well-defined.

A direct substitution into the defining relations
\eqref{eq:LL}--\eqref{eq:GG_same} shows that all brackets are preserved. For
example,
\[
\begin{aligned}
\{\sigma_\eta(G_r^+),\sigma_\eta(G_s^-)\}
&=\{G_{r-\eta}^+,G_{s+\eta}^-\} \\
&=-2L_{r+s}+(r-s-2\eta)H_{r+s} \\
&=\sigma_\eta\bigl(-2L_{r+s}+(r-s)H_{r+s}\bigr).
\end{aligned}
\]
The even and mixed relations are checked in the same way, using
\(\sigma_\eta(L_m)=L_m+\eta H_m\) and
\(\sigma_\eta(G_r^\pm)=G_{r\mp\eta}^\pm\). Since \(\sigma_{-\eta}\) is the
inverse of \(\sigma_\eta\), the map \(\sigma_\eta\) is a Lie superalgebra
isomorphism.
\end{proof}

By Proposition~\ref{prop:spectral-flow}, the Ramond and Neveu--Schwarz
algebras are isomorphic as Lie superalgebras. Since the realizations below are
naturally indexed by integers, from now on we work with the Ramond algebra
\(\mathfrak g=\mathfrak g[0]\).

Next we recall the weyl algebras and their simple modules.

\begin{defi}
The \emph{Weyl algebra} \(\mathcal D\) is the unital associative
\(\mathbb C\)-algebra generated by \(t^{\pm1}\) and \(D\), subject to
$
  tt^{-1}=t^{-1}t=1, [D,t]=t.
$
Equivalently, $
  Dt^m=t^m(D+m),m\in\mathbb Z,$
and every element of \(\mathcal D\) has a unique expression as a finite sum
$
  \sum_{m\in\mathbb Z} f_m(D)t^m,
 f_m(D)\in\mathbb C[D].
$
\end{defi}

A \(\mathcal D\)-module \(V\) is \(\mathbb C[t^{\pm1}]\)-\emph{torsion} if every
\(v\in V\) is annihilated by some nonzero element of \(\mathbb C[t^{\pm1}]\),
and \(\mathbb C[t^{\pm1}]\)-\emph{torsion-free} if no nonzero element of \(V\)
is annihilated by a nonzero element of \(\mathbb C[t^{\pm1}]\).

Let \(\mathbb C(t)[D]\) be the skew polynomial algebra generated by
\(\mathbb C(t)\) and \(D\), subject to
\[
  Df=fD+D(f)\quad(f\in\mathbb C(t)),
  \quad D(f)=t\frac{df}{dt}.
\]

\begin{theo}[{\cite{B,LZ}}]
\label{thm:D-modules}
Let \(V\) be a simple \(\mathcal D\)-module.  Then exactly one of the
following two cases occurs.

\begin{enumerate}
\item
\(V\) is \(\mathbb C[t^{\pm1}]\)-torsion-free.  In this case there exists a
non-unit irreducible element \(q\in\mathbb C(t)[D]\) such that
\[
  V\cong
  \mathcal D\big/
  \bigl(\mathcal D\cap \mathbb C(t)[D]q\bigr),
\]
where \(\mathbb C(t)[D]q\) denotes the left ideal generated by \(q\).

\item
\(V\) is \(\mathbb C[t^{\pm1}]\)-torsion.  In this case there exists
\(\lambda\in\mathbb C^*\) such that \(V\cong\Omega(\lambda)\), where
\(\Omega(\lambda)=\mathbb C[D]\) and
\[
  t^k\cdot f(D)=\lambda^k f(D-k),\quad
  D\cdot f(D)=Df(D)
\]
for all \(f(D)\in\mathbb C[D]\) and \(k\in\mathbb Z\).
\end{enumerate}
\end{theo}

The examples used later are the torsion modules \(\Omega(\lambda)\) and the
rank-one weight modules \(\mathcal V_a\); the latter are torsion-free as
\(\mathbb C[t^{\pm1}]\)-modules. We recall \(\mathcal V_a\) explicitly in
Section~\ref{subsec:rank-one-weight-modules}.

We now introduce the Weyl superalgebra with \(N\) odd variables.

\begin{defi}
Let
\[
  \mathcal A(1,N)=\mathbb C[t^{\pm1},\theta_1,\ldots,\theta_N]
\]
be the supercommutative algebra with one even invertible variable \(t\) and
odd variables \(\theta_1,\ldots,\theta_N\).  The \emph{Weyl superalgebra}
\(\mathcal{SD}(1,N)\) is the associative superalgebra generated by
\[
  t^{\pm1},\quad D=t\partial_t,\quad
  \theta_i,\quad \partial_i=\frac{\partial}{\partial\theta_i}
  \quad(1\le i\le N),
\]
where \(t^{\pm1}\) and \(D\) are even and \(\theta_i,\partial_i\) are odd,
subject to the relations
\[
  tt^{-1}=t^{-1}t=1,\ [D,t]=t,\
  [D,\theta_i]=[D,\partial_i]=[t^{\pm1},\theta_i]=[t^{\pm1},\partial_i]=0,
\]
and
\[
  \partial_i\theta_j+\theta_j\partial_i=\delta_{ij},\quad
  \theta_i\theta_j+\theta_j\theta_i=0,\quad
  \partial_i\partial_j+\partial_j\partial_i=0 .
\]
All odd derivatives are left odd derivatives.
\end{defi}

Let \(\mathcal C(N)\) be the finite-dimensional associative superalgebra
generated by \(\theta_i,\partial_i\) subject to the odd relations above. Then
\[
  \mathcal{SD}(1,N)\cong \mathcal D\otimes_{\mathbb C}\mathcal C(N),
\]
with \(\mathcal D\) purely even. In particular, \(\mathcal{SD}(1,N)\) has the
PBW-type basis
\[
  t^mD^r\theta_{i_1}\cdots\theta_{i_p}
  \partial_{j_1}\cdots\partial_{j_q},
\]
where \(m\in\mathbb Z\), \(r\in\mathbb N\),
\(1\le i_1<\cdots<i_p\le N\), and \(1\le j_1<\cdots<j_q\le N\).

Set \(\Lambda_N=\Lambda(\theta_1,\ldots,\theta_N)\). We regard \(\Lambda_N\) as
a standard \(\mathcal C(N)\)-module, with \(\theta_i\) acting by left
multiplication and \(\partial_i\) by left odd differentiation.

\begin{lemm}
\label{lem:Clifford-standard-module}
The natural representation
\[
  \rho:\mathcal C(N)\longrightarrow \End_{\mathbb C}(\Lambda_N)
\]
is an isomorphism of associative superalgebras.  Consequently,
\(\Lambda_N\) is, up to parity shift, the unique simple
\(\mathcal C(N)\)-module.
\end{lemm}
\begin{proof}
The monomials \(\theta_I\partial_J\), with
\(I,J\subseteq\{1,\ldots,N\}\), form a basis of \(\mathcal C(N)\), so
\[
  \dim \mathcal C(N)=2^{2N}=\dim \End_{\mathbb C}(\Lambda_N).
\]
It remains to check that the action on \(\Lambda_N\) is surjective. This is the
standard Fock representation of the finite Clifford algebra. The \(\theta_i\)
act by exterior multiplication and the \(\partial_i\) by contraction, and
together they generate all matrix units of \(\End_{\mathbb C}(\Lambda_N)\). For
instance, the projection onto the constant term is
\[
  P_\varnothing=\prod_{i=1}^N\partial_i\theta_i,
\]
and the operators \(\theta_I P_\varnothing\partial_J\), up to signs depending
on the chosen ordering, are the matrix units for the basis \(\{\theta_K\}\) of
\(\Lambda_N\). Thus \(\rho\) is surjective, and the dimension count makes it an
isomorphism.

The last assertion holds because a full matrix superalgebra has, up to parity
shift, a unique simple module, namely its natural module.
\end{proof}
\begin{prop}%[Simple \(\mathcal{SD}(1,N)\)-modules via Morita equivalence]
\label{pro:SD-modules_classification}
Let \(M\) be a simple \(\mathcal{SD}(1,N)\)-module.  Then, up to parity
shift, there exists a simple \(\mathcal D\)-module \(P\) such that
\(M\cong P\otimes\Lambda_N\) as \(\mathcal{SD}(1,N)\)-modules.  Conversely, if
\(P\) is a simple \(\mathcal D\)-module, then \(P\otimes\Lambda_N\) is a simple
\(\mathcal{SD}(1,N)\)-module.
\end{prop}
\begin{proof}
Put \(V=\Lambda_N\).  By Lemma~\ref{lem:Clifford-standard-module},
\[
  \mathcal C(N)\cong \End_{\mathbb C}(V),
  \quad
  \mathcal{SD}(1,N)\cong
  B:=\mathcal D\otimes_{\mathbb C}\End_{\mathbb C}(V).
\]
Choose a homogeneous basis \(v_1,\ldots,v_s\) of \(V\), where \(s=2^N\),
with corresponding matrix units \(E_{ij}\).

The argument below is the standard Morita equivalence between \(\mathcal D\)
and the matrix superalgebra \(B\); see, e.g., \cite{Lam}.
Let \(P\) be a simple \(\mathcal D\)-module. We first show that \(P\otimes V\)
is simple over \(B\). Let \(0\ne W\subseteq P\otimes V\) be a \(B\)-submodule,
and choose \(0\ne w=\sum_i p_i\otimes v_i\in W\) with some \(p_i\ne0\). Applying
\(1\otimes E_{ji}\) gives \(p_i\otimes v_j\in W\) for all \(j\). Since
\(\mathcal Dp_i=P\), it follows that \(P\otimes v_j\subseteq W\) for all \(j\),
and hence \(W=P\otimes V\). Thus \(P\otimes V\) is simple.

Conversely, let \(M\) be a simple \(B\)-module and set
\(e=1\otimes E_{11}\), so that \(eBe\cong\mathcal D\). Since
\(1=\sum_i1\otimes E_{ii}\), some \((1\otimes E_{ii})M\) is nonzero. Applying
\(1\otimes E_{1i}\) then shows that \(eM\ne0\). Put \(P=eM\). Then \(P\) is a
\(\mathcal D\)-module, and it is simple: if \(0\ne Q\subseteq P\) is a
\(\mathcal D\)-submodule, then \(BQ\) is a nonzero \(B\)-submodule of \(M\), so
\(BQ=M\), and applying \(e\) gives
\[
  P=eM=eBQ=eBe\,Q=Q.
\]

Now define
\[
  \Psi:P\otimes V\longrightarrow M,
  \quad
  \Psi(p\otimes v_i)=(1\otimes E_{i1})p .
\]
The matrix-unit relations show that \(\Psi\) is \(B\)-linear. Its image is
\(BP\), which is nonzero and hence equals \(M\). Suppose
\(\sum_i p_i\otimes v_i\in\ker\Psi\), so that
\(\sum_i (1\otimes E_{i1})p_i=0\). Applying \(1\otimes E_{1j}\) gives
\(p_j=0\) for every \(j\). Hence \(\Psi\) is injective, and \(M\cong P\otimes
V\).

Thus, up to the usual parity shift in the category of supermodules, every
simple \(\mathcal{SD}(1,N)\)-module is of the form \(P\otimes\Lambda_N\) for a
simple \(\mathcal D\)-module \(P\).
\end{proof}

\begin{coro}
\label{cor:standard-SD-functor}
For every \(N\ge0\), the assignment
\(\mathscr F_N(P)=P\otimes_{\mathbb C}\Lambda_N\) and
\(\mathscr F_N(f)=f\otimes\id_{\Lambda_N}\)
defines a functor
\[
  \mathscr F_N:
  \mathcal D\text{-}\mathrm{Mod}
  \longrightarrow
  \mathcal{SD}(1,N)\text{-}\mathrm{Mod}.
\]
Moreover, if \(P\) is simple as a \(\mathcal D\)-module, then
\(\mathscr F_N(P)\) is simple as a \(\mathcal{SD}(1,N)\)-module.
\end{coro}

\begin{proof}
The action is induced by
\(\mathcal{SD}(1,N)\cong\mathcal D\otimes_{\mathbb C}\mathcal C(N)\),
where \(\mathcal D\) acts on \(P\) and \(\mathcal C(N)\) on \(\Lambda_N\). If
\(f:P\to Q\) is a \(\mathcal D\)-module homomorphism, then
\(f\otimes\id_{\Lambda_N}\) commutes with both tensor-factor actions and is
therefore an \(\mathcal{SD}(1,N)\)-module homomorphism. The functorial
identities are immediate, and the simplicity statement follows from
Proposition~\ref{pro:SD-modules_classification}.
\end{proof}

%======================================================================
\section{Differential operator realizations of
\texorpdfstring{\(\mathfrak g\)}{g}}
\label{sec:Weyl_reps}
%======================================================================

This section fixes the differential-operator realizations used below. For
comparison, we first recall the one-odd-variable realization of \cite{CDLP2}.
We then turn to the standard two-odd-variable contact realization and to a
twisted form of it, obtained by an explicit change of variables. The twisted
form is equivalent to the contact realization, but its formulas are better
adapted to the four summands \(P,\ P\theta_1,\ P\theta_2,\ P\theta_1\theta_2\).

Throughout this section, set
\[
  A=\mathbb C[t^{\pm1},\theta_1,\theta_2],
  \quad
  W(1,2)=\Der(A),
  \quad
  D=t\partial_t,
  \quad
  \partial_i=\frac{\partial}{\partial\theta_i}\quad(i=1,2).
\]
We write \(\Theta=\theta_1\theta_2\), and all odd derivatives are left odd
derivatives. Put
\[
  E=\theta_2\partial_2+\theta_1\partial_1,
  \qquad
  J=\theta_2\partial_2-\theta_1\partial_1,
\]
so that
\[
  E(\theta_i)=\theta_i,\quad E(\Theta)=2\Theta,
  \qquad
  J(\theta_1)=-\theta_1,\quad J(\theta_2)=\theta_2,\quad J(\Theta)=0.
\]
For a \(\mathcal D\)-module \(P\), we write
\[
  \mathscr F_2(P)=P\otimes\Lambda(\theta_1,\theta_2)
  =P\oplus P\theta_1\oplus P\theta_2\oplus P\Theta ,
\]
so that every element has a unique expression
\[
  w=v_0+v_1\theta_1+v_2\theta_2+v_{12}\Theta,
  \quad v_0,v_1,v_2,v_{12}\in P .
\]
With our convention on left odd derivatives,
\(\partial_1\Theta=\theta_2\) and \(\partial_2\Theta=-\theta_1\), and the
operators \(t^{\pm1}\) and \(D\) act only on the coefficient in \(P\).

%----------------------------------------------------------------------
\subsection{The topological-twist realization}
\label{subsec:topological-modules}
%----------------------------------------------------------------------

Let \(W(1,1)=\Der(\mathbb C[t^{\pm1},\theta])\), with \(D=t\partial_t\) and
\(\partial_\theta=\partial/\partial\theta\).

\begin{prop}[cf.~\cite{CDLP2,EY90}]
\label{prop:topological_iso}
The assignment \(\pi^{\mathrm t}:\mathfrak g\to W(1,1)\) defined by
\[
\begin{aligned}
  \pi^{\mathrm t}(L_m)&=t^m\left(D+\frac m2\,\theta\partial_\theta\right),
&
  \pi^{\mathrm t}(H_m)&=t^m\theta\partial_\theta,\\
  \pi^{\mathrm t}(G_m^+)&=2t^m\theta D,
&
  \pi^{\mathrm t}(G_m^-)&=-t^m\partial_\theta,
\end{aligned}
\quad m\in\mathbb Z,
\]
is a Lie superalgebra isomorphism.
\end{prop}
\begin{proof}
This is the standard topological-twist realization of the centerless \(N=2\)
Ramond algebra in the convention \([L_m,L_n]=(n-m)L_{m+n}\). With left odd
differentiation, the displayed operators satisfy
\eqref{eq:LL}--\eqref{eq:GG_same}, and the identities
\[
  t^mD=\pi^{\mathrm t}(L_m)-\tfrac m2\pi^{\mathrm t}(H_m),\quad
  t^m\theta\partial_\theta=\pi^{\mathrm t}(H_m),
\]
\[
  t^m\theta D=\tfrac12\pi^{\mathrm t}(G_m^+),\quad
  t^m\partial_\theta=-\pi^{\mathrm t}(G_m^-)
\]
show that the map is surjective, hence an isomorphism.
\end{proof}

Via the embedding \(W(1,1)\hookrightarrow W(1,2)\) induced by
\(\theta\mapsto\theta_1\), we also regard \(\pi^{\mathrm t}\) as a realization
of \(\mathfrak g\) on \(A\). Here \(\theta_2\) is a spectator variable.

%----------------------------------------------------------------------
\subsection{The contact realization}
\label{subsec:contact-realization}
%----------------------------------------------------------------------

\begin{prop}[cf.~\cite{KL}]
\label{prop:contact-realization}
The assignment \(\pi^{\mathrm c}:\mathfrak g\to W(1,2)\) defined by
\begin{align}
  \pi^{\mathrm c}(L_m)
    &=t^m\left(D+\frac m2E\right),
  \label{eq:con-L}\\
  \pi^{\mathrm c}(H_m)
    &=t^mJ,
  \label{eq:con-H}\\
  \pi^{\mathrm c}(G_m^+)
    &=t^m\theta_2D-t^m\partial_1-mt^m\Theta\partial_1,
  \label{eq:con-Gp}\\
  \pi^{\mathrm c}(G_m^-)
    &=t^m\theta_1D-t^m\partial_2+mt^m\Theta\partial_2,
  \label{eq:con-Gm}
\end{align}
is a faithful Lie superalgebra homomorphism.
\end{prop}
\begin{proof}
These are the standard contact vector fields for the centerless \(N=2\) Ramond
algebra in the present sign convention. Using \([D,t^m]=mt^m\),
\(\partial_1\Theta=\theta_2\), \(\partial_2\Theta=-\theta_1\), and the
definitions of \(E\) and \(J\), one checks directly that the operators satisfy
\eqref{eq:LL}--\eqref{eq:GG_same}; see also \cite{KL}.

It remains to check faithfulness. A finite relation among the images splits by
parity and by \(t\)-degree. Fix \(m\). If
\(a\pi^{\mathrm c}(L_m)+b\pi^{\mathrm c}(H_m)=0\), then applying it to \(t\)
gives \(a=0\), and applying it to \(\theta_1\) gives \(b=0\). Similarly, if
\(c\pi^{\mathrm c}(G_m^+)+d\pi^{\mathrm c}(G_m^-)=0\), then applying it to
\(\theta_1\) and \(\theta_2\) gives
\[
  -c\,t^m(1+m\Theta)=0,\quad
  -d\,t^m(1-m\Theta)=0.
\]
Since \(t^m\) and \(1\pm m\Theta\) are invertible, \(c=d=0\). Hence the images
are linearly independent, and \(\pi^{\mathrm c}\) is faithful.
\end{proof}

%----------------------------------------------------------------------
\subsection{A twisted form of the contact realization}
\label{subsec:realization-asymmetric}
%----------------------------------------------------------------------

We now rewrite the contact realization by an explicit change of variables. The
resulting formulas are equivalent to those of the contact realization, but they
are better adapted to the four summands \(P,\ P\theta_1,\ P\theta_2,\ P\Theta\).
In particular, after the change of variables, \(G_m^+\) becomes the simple
operator \(t^m\partial_1\), while \(G_m^-\) carries the Weyl operator \(D\).

Every algebra automorphism \(\Phi\in\Aut_{\mathbb C}(A)\) induces, by
conjugation, an automorphism of the Weyl superalgebra,
\[
  T\longmapsto \Phi T\Phi^{-1},
  \quad T\in\mathcal{SD}(1,2).
\]
Hence, if \(\rho:\mathfrak g\to W(1,2)\) is a realization, then so is
\(X\mapsto \Phi\rho(X)\Phi^{-1}\).

\begin{prop}
\label{prop:contact-asym-conjugacy}
Define \(\pi^{\mathrm{tc}}:\mathfrak g\to W(1,2)\) by
\begin{subequations}\label{eq:realization-asym-formulas}
\begin{align}
  \pi^{\mathrm{tc}}(L_m)
    &=t^m\left(D+\frac m2E\right),
  \label{eq:asym-L}\\
  \pi^{\mathrm{tc}}(H_m)
    &=t^mJ,
  \label{eq:asym-H}\\
  \pi^{\mathrm{tc}}(G_m^+)
    &=t^m\partial_1,
  \label{eq:asym-Gp}\\
  \pi^{\mathrm{tc}}(G_m^-)
    &=-2t^m\theta_1D-t^m\partial_2-2mt^m\Theta\partial_2 .
  \label{eq:asym-Gm}
\end{align}
\end{subequations}
Let \(\Phi:A\to A\) be the algebra automorphism
\[
  \Phi(t)=t(1+\Theta),
  \quad
  \Phi(\theta_1)=-\theta_1,
  \quad
  \Phi(\theta_2)=\theta_2 .
\]
Then \(\Phi^2=\id_A\), and
\[
  \pi^{\mathrm c}(X)=\Phi\pi^{\mathrm{tc}}(X)\Phi^{-1},
  \quad X\in\mathfrak g .
\]
Consequently, \(\pi^{\mathrm{tc}}\) is a faithful Lie superalgebra
homomorphism. Moreover, the same linear map \(\Phi\) intertwines the two
\(\mathfrak g\)-actions on \(A\), so that
\(\Phi:A^{\mathrm{tc}}\to A^{\mathrm c}\) is a \(\mathfrak g\)-module
isomorphism.
\end{prop}

\begin{proof}
Since \(\Theta^2=0\) and \(\Phi(\Theta)=-\Theta\), we have
\[
  \Phi^2(t)=t(1+\Theta)(1-\Theta)=t,
  \quad
  \Phi^2(\theta_i)=\theta_i\quad(i=1,2),
\]
so \(\Phi^{-1}=\Phi\). Also
\[
  \Phi(t^m)=t^m(1+m\Theta),
  \quad
  \Phi(t^m\Theta)=-t^m\Theta .
\]

A direct calculation on the generators \(t,\theta_1,\theta_2\) gives
\[
\begin{aligned}
  \Phi D\Phi^{-1} &=D,
&
  \Phi E\Phi^{-1} &=E-2\Theta D,
&
  \Phi J\Phi^{-1} &=J,\\
  \Phi\partial_1\Phi^{-1} &=\theta_2D-\partial_1,
&
  \Phi\partial_2\Phi^{-1} &=\theta_1D+\partial_2 .
\end{aligned}
\]
For example,
\[
  \Phi\partial_1\Phi^{-1}(t)=t\theta_2,\quad
  \Phi\partial_1\Phi^{-1}(\theta_1)=-1,\quad
  \Phi\partial_1\Phi^{-1}(\theta_2)=0,
\]
which is the action of \(\theta_2D-\partial_1\). Similarly,
\[
  \Phi\partial_2\Phi^{-1}(t)=t\theta_1,\quad
  \Phi\partial_2\Phi^{-1}(\theta_1)=0,\quad
  \Phi\partial_2\Phi^{-1}(\theta_2)=1,
\]
which is the action of \(\theta_1D+\partial_2\); here the sign uses
\(\partial_2(\Theta)=-\theta_1\).

We now apply these identities to the generators. The images of \(E\) and \(J\)
lie in the ideal \((\theta_1,\theta_2)\), so multiplication by \(\Theta\)
annihilates both. Hence
\[
\begin{aligned}
  \Phi\pi^{\mathrm{tc}}(L_m)\Phi^{-1}
  &=
  t^m(1+m\Theta)
  \left(D+\frac m2(E-2\Theta D)\right)  \\
  &=
  t^m(1+m\Theta)
  \left(D+\frac m2E-m\Theta D\right)  \\
  &=
  t^m\left(D+\frac m2E\right)
  =
  \pi^{\mathrm c}(L_m),
\end{aligned}
\]
and similarly
\[
  \Phi\pi^{\mathrm{tc}}(H_m)\Phi^{-1}
  =
  t^m(1+m\Theta)J
  =
  t^mJ
  =
  \pi^{\mathrm c}(H_m).
\]

For \(G_m^+\),
\[
\begin{aligned}
  \Phi\pi^{\mathrm{tc}}(G_m^+)\Phi^{-1}
  &=
  t^m(1+m\Theta)(\theta_2D-\partial_1) \\
  &=
  t^m\theta_2D-t^m\partial_1-mt^m\Theta\partial_1
  =
  \pi^{\mathrm c}(G_m^+).
\end{aligned}
\]
Finally, using \(\Phi(t^m\theta_1)=-t^m\theta_1\) and
\(\Phi(t^m\Theta)=-t^m\Theta\), we obtain
\[
\begin{aligned}
  \Phi\pi^{\mathrm{tc}}(G_m^-)\Phi^{-1}
  &=
  2t^m\theta_1D
  -t^m(1+m\Theta)(\theta_1D+\partial_2)
  +2mt^m\Theta(\theta_1D+\partial_2)  \\
  &=
  t^m\theta_1D-t^m\partial_2+mt^m\Theta\partial_2
  =
  \pi^{\mathrm c}(G_m^-).
\end{aligned}
\]
Thus \(\pi^{\mathrm c}(X)=\Phi\pi^{\mathrm{tc}}(X)\Phi^{-1}\) for all
\(X\in\mathfrak g\). Since \(\pi^{\mathrm c}\) is faithful, so is
\(\pi^{\mathrm{tc}}\). Finally,
\[
  \Phi\pi^{\mathrm{tc}}(X)=\pi^{\mathrm c}(X)\Phi
  \quad(X\in\mathfrak g),
\]
so \(\Phi:A^{\mathrm{tc}}\to A^{\mathrm c}\) is a \(\mathfrak g\)-module
isomorphism.
\end{proof}

Let \(A^{\mathrm t}\) be the \(\mathfrak g\)-module obtained from the
topological-twist realization after the embedding \(\theta\mapsto\theta_1\),
and let \(A^{\mathrm{tc}}\) be the \(\mathfrak g\)-module obtained from
\(\pi^{\mathrm{tc}}\). Since \(A^{\mathrm{tc}}\cong A^{\mathrm c}\), it is
enough to compare \(A^{\mathrm t}\) with \(A^{\mathrm{tc}}\). We use the space
of invariants
\[
  {\rm H}^0(\mathfrak g,M)
  =\{\,u\in M\mid Xu=0\ \text{for all }X\in\mathfrak g\,\}.
\]

\begin{prop}
\label{prop:topological-comparison}
The topological-twist module \(A^{\mathrm t}\) is not isomorphic to the
contact module \(A^{\mathrm c}\).
\end{prop}

\begin{proof}
Since \({\rm H}^0(\mathfrak g,-)\) is functorial, the dimension of
\({\rm H}^0(\mathfrak g,M)\) is an isomorphism invariant.

For the topological-twist action on \(A\), we have
\(G_0^-=-\partial_1\) and \(G_0^+=2\theta_1D\). Let \(f\in A^{\mathrm t}\) be
invariant. The equation \(G_0^-f=0\) gives \(\partial_1f=0\). Since
\(\theta_2\) is a spectator in the embedded topological-twist realization, this
means \(f=a(t)+b(t)\theta_2\). Then \(G_0^+f=0\) gives \(Da(t)=0\) and
\(Db(t)=0\), so \(a(t),b(t)\in\mathbb C\), and
\[
  {\rm H}^0(\mathfrak g,A^{\mathrm t})
  =
  \mathbb C1\oplus\mathbb C\theta_2 .
\]

For the action \(\pi^{\mathrm{tc}}\), write
\(f=a(t)+b(t)\theta_1+c(t)\theta_2+d(t)\Theta\). Since \(G_0^+=\partial_1\),
invariance gives \(b(t)=d(t)=0\), so \(f=a(t)+c(t)\theta_2\). Since
\(H_0=J\), the equation \(H_0f=0\) gives \(c(t)=0\). Finally \(L_0=D\) gives
\(Da(t)=0\), hence \(a(t)\in\mathbb C\). Therefore
\[
  {\rm H}^0(\mathfrak g,A^{\mathrm{tc}})=\mathbb C1 .
\]
The two invariant spaces have different dimensions, so
\(A^{\mathrm t}\not\cong A^{\mathrm{tc}}\). Since
\(A^{\mathrm{tc}}\cong A^{\mathrm c}\), we also have
\(A^{\mathrm t}\not\cong A^{\mathrm c}\).
\end{proof}

%======================================================================
\section{\texorpdfstring{\(A\)-valued \(1\)-cocycles and a
two-parameter family of \(\mathfrak g\)-modules}
{A-valued 1-cocycles and a two-parameter family of g-modules}}
\label{sec:twisting-functors}
%======================================================================

This section constructs the modules studied in the rest of the paper. We first
compute the \(A\)-valued \(1\)-cocycles on \(\mathfrak g\), then use them to
deform the realization \(\pi^{\mathrm{tc}}\) into a two-parameter family of
\(\mathfrak g\)-modules \(F_{\alpha,\beta}\).

We pass from the realization \(\pi^{\mathrm{tc}}\) to the \(\mathfrak g\)-modules
studied below. Put \(A=\mathbb C[t^{\pm1},\theta_1,\theta_2]\) and
\(\Theta=\theta_1\theta_2\). The realization
\(\pi^{\mathrm{tc}}:\mathfrak g\to\Der(A)\) makes \(A\) a \(\mathfrak g\)-module,
denoted \(A^{\mathrm{tc}}\). For a \(\mathcal D\)-module \(P\), the space
\[
  \mathscr F_2(P)=P\otimes\Lambda(\theta_1,\theta_2)
\]
carries the induced \(\mathfrak g\)-action through the embedding
\(\pi^{\mathrm{tc}}(\mathfrak g)\subset \mathcal{SD}(1,2)\), and we deform this
action by \(A\)-valued \(1\)-cocycles. Throughout this section, elements of
\(A\) act on \(\mathscr F_2(P)\) through the natural \(A\)-action: the Laurent
variable \(t\) acts through the \(\mathcal D\)-module structure on \(P\), and
the odd variables act on the exterior factor.

%----------------------------------------------------------------------
\subsection{\texorpdfstring{\(A\)-valued \(1\)-cocycles}
{A-valued 1-cocycles}}
%----------------------------------------------------------------------

\begin{lemm}
\label{lem:cocycle-perturbation}
Let \(d:\mathfrak g\to A^{\mathrm{tc}}\) be an even \(1\)-cocycle.  For
every \(\mathcal D\)-module \(P\), define
\[
  X\cdot_d w=\pi^{\mathrm{tc}}(X)w+d(X)w,
  \quad X\in\mathfrak g,\quad w\in\mathscr F_2(P).
\]
Then \(\cdot_d\) is a \(\mathfrak g\)-module structure on
\(\mathscr F_2(P)\).
\end{lemm}

\begin{proof}
Let \(X,Y\in\mathfrak g\) be homogeneous. For homogeneous \(a\in A\),
multiplication by \(a\) satisfies
\[
  [\pi^{\mathrm{tc}}(X),a]=\pi^{\mathrm{tc}}(X)(a),
\]
the right-hand side again being viewed as a multiplication operator. Since
\(d\) is even, \(d(X)\) has parity \(|X|\). As \(A\) is supercommutative, the
multiplication operators \(d(X)\) and \(d(Y)\) supercommute. Therefore
\[
\begin{aligned}
&[\pi^{\mathrm{tc}}(X)+d(X),\,\pi^{\mathrm{tc}}(Y)+d(Y)]  \\
&\quad =
  \pi^{\mathrm{tc}}([X,Y])
  +\pi^{\mathrm{tc}}(X)d(Y)
  -(-1)^{|X||Y|}\pi^{\mathrm{tc}}(Y)d(X),
\end{aligned}
\]
and the last two terms equal \(d([X,Y])\) by the cocycle identity. Hence the
operators \(X\mapsto \pi^{\mathrm{tc}}(X)+d(X)\) satisfy the defining bracket
relations of \(\mathfrak g\).
\end{proof}

\begin{lemm}
\label{lem:coboundary-weyl-twist}
Let \(d:\mathfrak g\to A^{\mathrm{tc}}\) be an even \(1\)-cocycle and let
\(a\in A_{\bar0}\).  Regard \(a\) as a multiplication operator on \(A\), and
put
\[
  \tau_a=\exp(-\operatorname{ad}a)\in\Aut\mathcal{SD}(1,2).
\]
Then, for every \(X\in\mathfrak g\),
\[
  \tau_a\bigl(\pi^{\mathrm{tc}}(X)+d(X)\bigr)
  =
  \pi^{\mathrm{tc}}(X)+d(X)+(\delta a)(X),
  \quad
  (\delta a)(X)=\pi^{\mathrm{tc}}(X)a.
\]
Consequently, for every \(\mathcal{SD}(1,2)\)-module \(M\), the
\(\mathfrak g\)-module defined by the cocycle \(d+\delta a\) is the
restriction of the twisted Weyl-superalgebra module \(M^{\tau_a}\) along
the homomorphism \(X\mapsto \pi^{\mathrm{tc}}(X)+d(X)\).
\end{lemm}

\begin{proof}
Since \(a\) is a multiplication operator, \(\operatorname{ad}a\) lowers the
order of differential operators. Hence it is locally nilpotent on
\(\mathcal{SD}(1,2)\), and \(\tau_a\) is well-defined. For every derivation
\(T\in W(1,2)\) one has \([a,T]=-T(a)\), with \(T(a)\) viewed as a
multiplication operator, so
\[
  \tau_a(T)=\exp(-\operatorname{ad}a)(T)=T+T(a).
\]
Taking \(T=\pi^{\mathrm{tc}}(X)\) gives
\[
  \tau_a(\pi^{\mathrm{tc}}(X))
  =
  \pi^{\mathrm{tc}}(X)+\pi^{\mathrm{tc}}(X)a
  =
  \pi^{\mathrm{tc}}(X)+(\delta a)(X).
\]
Since \(\tau_a\) fixes multiplication operators from \(A\), we have
\(\tau_a(d(X))=d(X)\), and the formula follows. The final assertion is this
identity applied to the action operators on \(M^{\tau_a}\).
\end{proof}

We now compute the \(A\)-valued \(1\)-cocycles needed below.

\begin{prop}
\label{prop:multiplication-cocycles}
Every odd \(1\)-cocycle \(d:\mathfrak g\to A^{\mathrm{tc}}\) is a
coboundary.  Every even \(1\)-cocycle \(d:\mathfrak g\to A^{\mathrm{tc}}\)
has the form
\[
  d=d_{\alpha,\beta,\gamma}+\delta a,
  \quad
  \alpha,\beta,\gamma\in\mathbb C,\quad a\in A_{\bar0},
\]
where \((\delta a)(X)=\pi^{\mathrm{tc}}(X)a\) and
\begin{subequations}\label{eq:cocycle-alpha-beta-gamma}
\begin{align}
  d_{\alpha,\beta,\gamma}(L_m)
    &=(\gamma+m\alpha)t^m
      +\frac{m^2}{2}(2\alpha-\beta)t^m\Theta,
      \label{eq:cocycle-abg-L}\\
  d_{\alpha,\beta,\gamma}(H_m)
    &=\beta t^m-m(2\alpha-\beta)t^m\Theta,
      \label{eq:cocycle-abg-H}\\
  d_{\alpha,\beta,\gamma}(G_m^+)
    &=m(2\alpha-\beta)t^m\theta_2,
      \label{eq:cocycle-abg-Gp}\\
  d_{\alpha,\beta,\gamma}(G_m^-)
    &=\bigl(-2\gamma-m(2\alpha+\beta)\bigr)t^m\theta_1 .
      \label{eq:cocycle-abg-Gm}
\end{align}
\end{subequations}
For fixed \(\alpha,\beta,\gamma\), the element \(a\) is unique modulo
\(\mathbb C1\).  Consequently,
\[
  {\rm H}^1(\mathfrak g,A^{\mathrm{tc}})_{\bar0}\cong\mathbb C^3,
  \quad
  {\rm H}^1(\mathfrak g,A^{\mathrm{tc}})_{\bar1}=0.
\]
\end{prop}

\begin{proof}
We use the \(\mathbb Z\times\mathbb Z\)-grading given by the
\(L_0\)-degree and the \(H_0\)-charge.  On \(A\),
\[
  \pi^{\mathrm{tc}}(L_0)=D,
  \quad
  \pi^{\mathrm{tc}}(H_0)=J,
\]
so \(t^p\) has \(L_0\)-degree \(p\), while
\[
  \operatorname{ch}\theta_1=-1,\quad
  \operatorname{ch}\theta_2=1,\quad
  \operatorname{ch}\Theta=0.
\]
The generators of \(\mathfrak g\) have \(L_0\)-degree \(m\) and charges
\[
  \operatorname{ch}L_m=\operatorname{ch}H_m=0,\quad
  \operatorname{ch}G_m^+=1,\quad
  \operatorname{ch}G_m^-=-1 .
\]
Since \(L_0\) and \(H_0\) act semisimply on both \(\mathfrak g\) and \(A\),
the cocycle complex decomposes into homogeneous bidegree parts. It therefore
suffices to compute homogeneous cocycles.

Let \(d\) have bidegree \((p,s)\). If \(X\) has \(L_0\)-degree \(q\), the
cocycle identity for \((L_0,X)\) gives
\[
  qd(X)=(q+p)d(X)-\pi^{\mathrm{tc}}(X)d(L_0),
\]
so \(d=\delta(p^{-1}d(L_0))\) whenever \(p\ne0\). Similarly, the pair
\((H_0,X)\) gives
\[
  rd(X)=(r+s)d(X)-\pi^{\mathrm{tc}}(X)d(H_0)
\]
for \(X\) of \(H_0\)-charge \(r\), so \(d=\delta(s^{-1}d(H_0))\) whenever
\(p=0\) and \(s\ne0\). Therefore only bidegree \((0,0)\) can contribute to
cohomology.

Consider next the odd cocycles of bidegree \((0,0)\). The charge condition
forces \(d(L_m)=d(H_m)=0\), since there are no odd monomials of degree \(m\) and
charge \(0\), and \(d(G_m^\pm)=0\), since there are no even monomials of degree
\(m\) and charge \(\pm1\). Thus every odd bidegree-zero cocycle is zero, and
every odd cocycle is a coboundary.

It remains to compute the even bidegree-zero cocycles. Write
\[
\begin{aligned}
  d(L_m)&=a_m t^m+b_m t^m\Theta,
&
  d(H_m)&=c_m t^m+e_m t^m\Theta,\\
  d(G_m^+)&=v_m t^m\theta_2,
&
  d(G_m^-)&=r_m t^m\theta_1 .
\end{aligned}
\]
Applying the cocycle identity to
\([L_m,L_n]=(n-m)L_{m+n}\) and \([L_m,H_n]=nH_{m+n}\) gives
\[
\begin{gathered}
  (n-m)a_{m+n}=n a_n-m a_m,
  \quad
  (n-m)b_{m+n}=(m+n)(b_n-b_m),\\
  n c_{m+n}=n c_n,
  \quad
  n e_{m+n}=(m+n)e_n ,
\end{gathered}
\]
so that
\[
  a_m=\gamma+\alpha m,\quad
  c_m=\beta,\quad
  e_m=\kappa m
\]
for some \(\alpha,\beta,\gamma,\kappa\in\mathbb C\).

The relations \([H_m,G_n^+]=G_{m+n}^+\) and \([H_m,G_n^-]=-G_{m+n}^-\) yield
\[
  v_{m+n}=v_n-e_m,
  \quad
  r_{m+n}=r_n-2mc_m+e_m,
\]
hence
\[
  v_m=\zeta-\kappa m,
  \quad
  r_m=\eta+m(\kappa-2\beta)
\]
for some \(\zeta,\eta\in\mathbb C\). The relation
\([L_m,G_n^+]=\left(n-\frac m2\right)G_{m+n}^+\) gives
\[
  \left(n-\frac m2\right)v_{m+n}
  =
  \left(n+\frac m2\right)v_n-b_m,
\]
and therefore
\[
  b_m=\zeta m-\frac{\kappa}{2}m^2 .
\]
Finally, the scalar component of
\(\{G_m^+,G_n^-\}=-2L_{m+n}+(m-n)H_{m+n}\) gives
\[
  r_n-v_m=-2a_{m+n}+(m-n)c_{m+n},
\]
and substituting yields
\[
  \eta=\zeta-2\gamma,
  \quad
  \kappa=\beta-2\alpha .
\]
Thus
\[
\begin{gathered}
  e_m=-m(2\alpha-\beta),
  \quad
  v_m=\zeta+m(2\alpha-\beta),\\
  b_m=\zeta m+\frac{m^2}{2}(2\alpha-\beta),
  \quad
  r_m=\zeta-2\gamma-m(2\alpha+\beta),
\end{gathered}
\]
and the remaining defining relations impose no further conditions.

The parameter \(\zeta\) is a coboundary. Indeed,
\[
\begin{gathered}
  (\delta\Theta)(L_m)=m t^m\Theta,
  \quad
  (\delta\Theta)(H_m)=0,\\
  (\delta\Theta)(G_m^+)=t^m\theta_2,
  \quad
  (\delta\Theta)(G_m^-)=t^m\theta_1 ,
\end{gathered}
\]
so subtracting \(\zeta\delta\Theta\) removes exactly the \(\zeta\)-terms and
leaves the representative \(d_{\alpha,\beta,\gamma}\) above.

It remains to check that the three displayed representatives are not
coboundaries. A bidegree-zero even coboundary can come only from
\[
  A_{\bar0}^{(0,0)}=\mathbb C1\oplus\mathbb C\Theta .
\]
The constant gives zero, and \(\Theta\) gives the coboundary already used.
Hence no nonzero linear combination of \(d_{1,0,0}\), \(d_{0,1,0}\),
\(d_{0,0,1}\) is a coboundary, and
\[
  {\rm H}^1(\mathfrak g,A^{\mathrm{tc}})_{\bar0}\cong\mathbb C^3,
  \quad
  {\rm H}^1(\mathfrak g,A^{\mathrm{tc}})_{\bar1}=0.
\]

Finally, \(\delta a=0\) if and only if
\(a\in{\rm H}^0(\mathfrak g,A^{\mathrm{tc}})\). By
Proposition~\ref{prop:topological-comparison}, this invariant space is
\(\mathbb C1\). Hence, once the representative \(d_{\alpha,\beta,\gamma}\) is
fixed, the coboundary potential \(a\) is unique modulo constants.
\end{proof}

We shall write \(d_{\alpha,\beta}:=d_{\alpha,\beta,0}\), so that
\begin{subequations}\label{eq:cocycle-alpha-beta}
\begin{align}
  d_{\alpha,\beta}(L_m)
    &=m\alpha\,t^m
      +\frac{m^2}{2}(2\alpha-\beta)t^m\Theta,
      \label{eq:cocycle-alpha-beta-L}\\
  d_{\alpha,\beta}(H_m)
    &=\beta t^m-m(2\alpha-\beta)t^m\Theta,
      \label{eq:cocycle-alpha-beta-H}\\
  d_{\alpha,\beta}(G_m^+)
    &=m(2\alpha-\beta)t^m\theta_2,
      \label{eq:cocycle-alpha-beta-Gp}\\
  d_{\alpha,\beta}(G_m^-)
    &=-m(2\alpha+\beta)t^m\theta_1 .
      \label{eq:cocycle-alpha-beta-Gm}
\end{align}
\end{subequations}

%----------------------------------------------------------------------
\subsection{Reduction to two parameters}
\label{subsec:reduction-two-parameters}
%----------------------------------------------------------------------

The preceding proposition gives three normalized even parameters. We now
explain why only two of them need to be kept in the module construction.

Let \(\partial_D\) be the linear map on the generators of \(\mathcal{SD}(1,2)\)
defined by
\[
  \partial_D(D)=1,
  \quad
  \partial_D(t^{\pm1})=0,
  \quad
  \partial_D(\theta_i)=\partial_D(\partial_i)=0
  \quad(i=1,2).
\]
This assignment preserves the defining relations of \(\mathcal{SD}(1,2)\)---for
example, \(\partial_D([D,t]-t)=[1,t]=0\)---so \(\partial_D\) is a derivation of
\(\mathcal{SD}(1,2)\). Since \(\partial_D\) lowers the degree in \(D\), it is
locally nilpotent, and
\[
  \sigma_\gamma:=\exp(\gamma\partial_D)
\]
is an automorphism of \(\mathcal{SD}(1,2)\). Explicitly,
\[
  \sigma_\gamma(D)=D+\gamma,\quad
  \sigma_\gamma(t^{\pm1})=t^{\pm1},\quad
  \sigma_\gamma(\theta_i)=\theta_i,\quad
  \sigma_\gamma(\partial_i)=\partial_i .
\]
For \(a\in A_{\bar0}\), let \(\tau_a\) be the automorphism of
\(\mathcal{SD}(1,2)\) from Lemma~\ref{lem:coboundary-weyl-twist}.

\begin{prop}
\label{prop:gamma-coboundary-automorphisms}
For every \(X\in\mathfrak g\),
\[
  \sigma_\gamma(\pi^{\mathrm{tc}}(X))
  =
  \pi^{\mathrm{tc}}(X)+\gamma d_{0,0,1}(X),
  \quad
  \tau_a(\pi^{\mathrm{tc}}(X))
  =
  \pi^{\mathrm{tc}}(X)+(\delta a)(X),
\]
and
\[
  \tau_a\sigma_\gamma
  \bigl(\pi^{\mathrm{tc}}(X)+d_{\alpha,\beta}(X)\bigr)
  =
  \pi^{\mathrm{tc}}(X)+d_{\alpha,\beta,\gamma}(X)+(\delta a)(X).
\]
Thus every cocycle of the form \(d_{\alpha,\beta,\gamma}+\delta a\) is
obtained from \(d_{\alpha,\beta}\) by applying an automorphism of
\(\mathcal{SD}(1,2)\).
\end{prop}

\begin{proof}
The first identity follows directly from the formulas for
\(\pi^{\mathrm{tc}}\). Since only \(D\) is changed,
\[
  \sigma_\gamma(\pi^{\mathrm{tc}}(L_m))
  =
  t^m\left(D+\gamma+\frac m2E\right)
  =
  \pi^{\mathrm{tc}}(L_m)+\gamma t^m,
\]
and
\[
\begin{aligned}
  \sigma_\gamma(\pi^{\mathrm{tc}}(G_m^-))
  &=
  -2t^m\theta_1(D+\gamma)
  -t^m\partial_2
  -2mt^m\Theta\partial_2        \\
  &=
  \pi^{\mathrm{tc}}(G_m^-)-2\gamma t^m\theta_1 ,
\end{aligned}
\]
while the images of \(H_m\) and \(G_m^+\) are unchanged. These four identities
are exactly
\(\sigma_\gamma(\pi^{\mathrm{tc}}(X))=\pi^{\mathrm{tc}}(X)+\gamma d_{0,0,1}(X)\).
The second identity is Lemma~\ref{lem:coboundary-weyl-twist} applied to the
zero cocycle. Since both \(\sigma_\gamma\) and \(\tau_a\) fix multiplication
operators from \(A\), we have
\(\tau_a\sigma_\gamma(d_{\alpha,\beta}(X))=d_{\alpha,\beta}(X)\), and combining
the identities gives the last formula.
\end{proof}

If \(S\) is an associative algebra, \(\varphi\in\Aut(S)\), and \(M\) is an
\(S\)-module, let \(M^\varphi\) denote the \(S\)-module with the same
underlying vector space and action \(s\cdot_\varphi v=\varphi(s)v\).

\begin{coro}
\label{cor:gamma-coboundary-module-reduction}
Let \(P\) be a \(\mathcal D\)-module.

\begin{enumerate}
\item
Let \(P^{\langle\gamma\rangle}\) be the \(\mathcal D\)-module with the same
\(t^{\pm1}\)-action as \(P\), but with \(D\) acting as \(D+\gamma\).  Then
there is a natural \(\mathcal{SD}(1,2)\)-module isomorphism
\[
  \mathscr F_2(P)^{\sigma_\gamma}
  \cong
  \mathscr F_2(P^{\langle\gamma\rangle}),
\]
where the two sides are identified as vector spaces.

\item
The \(\mathfrak g\)-module obtained from the cocycle
\(d_{\alpha,\beta,\gamma}\) on \(\mathscr F_2(P)\) is naturally identified
with the \(\mathfrak g\)-module obtained from the cocycle
\(d_{\alpha,\beta}\) on \(\mathscr F_2(P^{\langle\gamma\rangle})\).

\item
More generally, the \(\mathfrak g\)-module obtained from
\(d_{\alpha,\beta,\gamma}+\delta a\) on \(\mathscr F_2(P)\) is the
restriction of
\[
  \mathscr F_2(P)^{\tau_a\sigma_\gamma}
\]
along the two-parameter homomorphism
\(X\mapsto \pi^{\mathrm{tc}}(X)+d_{\alpha,\beta}(X)\).
\end{enumerate}

Consequently, the parameter \(\gamma\) is absorbed by shifting the
\(D\)-action of \(P\), while every coboundary term is absorbed by twisting
the underlying Weyl-superalgebra module. Thus, for the module construction
below, it is enough to keep the two normalized parameters \(\alpha\) and
\(\beta\).
\end{coro}

\begin{proof}
The automorphism \(\sigma_\gamma\) fixes \(t^{\pm1}\), \(\theta_i\), and
\(\partial_i\) and sends \(D\) to \(D+\gamma\), so twisting
\(\mathscr F_2(P)\) by \(\sigma_\gamma\) changes only the action of \(D\) on
the Laurent Weyl factor. This gives
\[
  \mathscr F_2(P)^{\sigma_\gamma}
  \cong
  \mathscr F_2(P^{\langle\gamma\rangle}).
\]
By Proposition~\ref{prop:gamma-coboundary-automorphisms},
\(\sigma_\gamma(\pi^{\mathrm{tc}}(X))=\pi^{\mathrm{tc}}(X)+\gamma
d_{0,0,1}(X)\), and since \(\sigma_\gamma\) fixes multiplication operators
from \(A\), it fixes \(d_{\alpha,\beta}(X)\). Therefore the action
\(X\mapsto \pi^{\mathrm{tc}}(X)+d_{\alpha,\beta,\gamma}(X)\) on
\(\mathscr F_2(P)\) is identified with the action
\(X\mapsto \pi^{\mathrm{tc}}(X)+d_{\alpha,\beta}(X)\) on
\(\mathscr F_2(P^{\langle\gamma\rangle})\). The last assertion follows from
the final identity in
Proposition~\ref{prop:gamma-coboundary-automorphisms}.
\end{proof}

\begin{prop}
\label{prop:F-alpha-beta-module}
Let \(\alpha,\beta\in\mathbb C\).  For every \(\mathcal D\)-module \(P\),
the following formulas define a \(\mathfrak g\)-module structure on
\(P\otimes\Lambda(\theta_1,\theta_2)\):
\begin{align}
  L_m
  &=
  t^m\left(D+\frac m2E\right)
  +m\alpha\,t^m
  +\frac{m^2}{2}(2\alpha-\beta)t^m\Theta,
  \label{eq:action-L-alpha-beta}\\
  H_m
  &=
  t^mJ+\beta t^m
  -m(2\alpha-\beta)t^m\Theta,
  \label{eq:action-H-alpha-beta}\\
  G_m^+
  &=
  t^m\partial_1
  +m(2\alpha-\beta)t^m\theta_2,
  \label{eq:action-Gp-alpha-beta}\\
  G_m^-
  &=
  -2t^m\theta_1D
  -t^m\partial_2
  -2m\,t^m\Theta\partial_2
  -m(2\alpha+\beta)t^m\theta_1 ,
  \label{eq:action-Gm-alpha-beta}
\end{align}
where
\[
  E=\theta_1\partial_1+\theta_2\partial_2,
  \quad
  J=\theta_2\partial_2-\theta_1\partial_1,
  \quad
  \Theta=\theta_1\theta_2.
\]
We denote the resulting module by \(F_{\alpha,\beta}(P)\).
\end{prop}

\begin{proof}
The displayed operators are obtained by adding the cocycle
\(d_{\alpha,\beta}\) to the action \(\pi^{\mathrm{tc}}\), so the claim
follows from Lemma~\ref{lem:cocycle-perturbation}.
\end{proof}

\begin{coro}
\label{cor:F-alpha-beta-functor}
For fixed \(\alpha,\beta\in\mathbb C\), the assignment
\(P\mapsto F_{\alpha,\beta}(P)\) defines a functor
\[
  F_{\alpha,\beta}:
  \mathcal D\text{-}\mathrm{Mod}
  \longrightarrow
  \mathfrak g\text{-}\mathrm{Mod},
\]
whose value on a \(\mathcal D\)-module homomorphism \(f:P\to Q\) is
\(f\otimes\id_{\Lambda(\theta_1,\theta_2)}\).
\end{coro}

\begin{proof}
Since \(f:P\to Q\) is a \(\mathcal D\)-module homomorphism, it commutes with
\(t^{\pm1}\) and \(D\), while the operators
\(\theta_i,\partial_i,E,J,\Theta\) act on the exterior factor and hence
commute with \(f\otimes\id_{\Lambda(\theta_1,\theta_2)}\). Therefore
\(f\otimes\id_{\Lambda(\theta_1,\theta_2)}\) intertwines the four operators in
\eqref{eq:action-L-alpha-beta}--\eqref{eq:action-Gm-alpha-beta} and so defines
a \(\mathfrak g\)-module homomorphism. The functorial identities are
immediate.
\end{proof}

%======================================================================
\section{Irreducibility and isomorphism classes}
\label{sec:irreducibility}
%======================================================================

This section analyzes the modules \(F_{\alpha,\beta}(P)\). We give a simplicity
criterion in terms of the two linear forms \(2\alpha\pm\beta\), classify the
simple modules up to isomorphism, and describe the structure on the two
singular hyperplanes.

We determine the irreducibility and the generic isomorphism classes of
\(F_{\alpha,\beta}(P)\). Unless otherwise stated, \(P\) is a simple
\(\mathcal D\)-module in the irreducibility and isomorphism results; the
singular exact sequences at the end are valid for arbitrary
\(\mathcal D\)-modules.

Throughout this section, set
\[
  M=F_{\alpha,\beta}(P)=P\otimes\Lambda,
  \quad
  \Lambda=\Lambda(\theta_1,\theta_2),
  \quad
  \Theta=\theta_1\theta_2 ,
\]
identify \(P\) with \(P\otimes1\), and put
\[
  \mu_-=2\alpha-\beta,
  \qquad
  \mu_+=2\alpha+\beta .
\]
The two linear forms \(\mu_-\) and \(\mu_+\) control the submodule structure.
We first record the component formulas used in the proofs.

\begin{lemm}
\label{lem:odd-component-formulas}
With respect to the decomposition
\[
  M=P\oplus P\theta_1\oplus P\theta_2\oplus P\Theta,
\]
the odd generators act as follows. For every \(v\in P\) and
\(m\in\mathbb Z\),
\[
\begin{aligned}
  G_m^+v&=m\mu_-t^mv\theta_2,
&
  G_m^+(v\theta_1)&=t^mv-m\mu_-t^mv\Theta,\\
  G_m^+(v\theta_2)&=0,
&
  G_m^+(v\Theta)&=t^mv\theta_2,
\end{aligned}
\]
and
\[
\begin{aligned}
  G_m^-v&=-t^m(2D+m\mu_+)v\theta_1,
&
  G_m^-(v\theta_1)&=0,\\
  G_m^-(v\theta_2)&=-t^mv-t^m(2D+2m+m\mu_+)v\Theta,
&
  G_m^-(v\Theta)&=t^mv\theta_1 .
\end{aligned}
\]
Here \(t^m\) and \(D\) act on the coefficient in \(P\). In particular,
\[
\begin{gathered}
  G_0^+v=0,\quad G_0^+(v\theta_1)=v,\quad
  G_0^+(v\theta_2)=0,\quad G_0^+(v\Theta)=v\theta_2,\\
  G_0^-v=-2Dv\theta_1,\quad G_0^-(v\theta_1)=0,\quad
  G_0^-(v\theta_2)=-v-2Dv\Theta,\quad
  G_0^-(v\Theta)=v\theta_1 .
\end{gathered}
\]
Moreover, the four exterior components are \(H_0\)-eigenspaces, with
\[
  H_0|_P=\beta,\quad
  H_0|_{P\theta_1}=\beta-1,\quad
  H_0|_{P\theta_2}=\beta+1,\quad
  H_0|_{P\Theta}=\beta .
\]
\end{lemm}

\begin{proof}
For \(G_m^+\), use
\[
  G_m^+=t^m\partial_1+m\mu_-t^m\theta_2
\]
together with
\[
  \partial_1(1)=0,\quad
  \partial_1\theta_1=1,\quad
  \partial_1\theta_2=0,\quad
  \partial_1\Theta=\theta_2,
\]
and \(\theta_2^2=0\), \(\theta_2\theta_1=-\Theta\). For \(G_m^-\), use
\[
  G_m^-=-2t^m\theta_1D-t^m\partial_2
         -2mt^m\Theta\partial_2-m\mu_+t^m\theta_1,
\]
together with
\[
  \partial_2(1)=\partial_2\theta_1=0,\quad
  \partial_2\theta_2=1,\quad
  \partial_2\Theta=-\theta_1,
\]
and \(\theta_1^2=0\), \(\theta_1\Theta=0\). The charge statement follows from
\(H_0=J+\beta\) and
\[
  J(1)=J(\Theta)=0,\quad
  J(\theta_1)=-\theta_1,\quad
  J(\theta_2)=\theta_2 .
\]
\end{proof}

\begin{lemm}
\label{lem:even-component-formulas}
With respect to the same decomposition, the even generators act as follows.
For every \(v\in P\) and \(m\in\mathbb Z\),
\[
\begin{aligned}
  L_mv&=t^m(D+m\alpha)v+\frac{m^2}{2}\mu_-t^mv\Theta,
&
  L_m(v\theta_1)&=t^m\bigl(D+m(\alpha+\tfrac12)\bigr)v\theta_1,\\
  L_m(v\theta_2)&=t^m\bigl(D+m(\alpha+\tfrac12)\bigr)v\theta_2,
&
  L_m(v\Theta)&=t^m(D+m(\alpha+1))v\Theta,
\end{aligned}
\]
and
\[
\begin{aligned}
  H_mv&=\beta t^mv-m\mu_-t^mv\Theta,
&
  H_m(v\theta_1)&=(\beta-1)t^mv\theta_1,\\
  H_m(v\theta_2)&=(\beta+1)t^mv\theta_2,
&
  H_m(v\Theta)&=\beta t^mv\Theta .
\end{aligned}
\]
\end{lemm}

\begin{proof}
The formulas follow from
\[
  L_m=t^m\left(D+\frac m2E\right)+m\alpha t^m
      +\frac{m^2}{2}\mu_-t^m\Theta,
  \quad
  H_m=t^mJ+\beta t^m-m\mu_-t^m\Theta,
\]
together with the exterior eigenvalues
\[
  E=(0,1,1,2),
  \quad
  J=(0,-1,1,0)
\]
on \(1,\theta_1,\theta_2,\Theta\), respectively. Since
\(\Theta\cdot1=\Theta\) and
\(\Theta\theta_1=\Theta\theta_2=\Theta^2=0\), substitution gives the result.
\end{proof}

%----------------------------------------------------------------------
\subsection{Irreducibility}
\label{subsec:generic_irreducibility}
%----------------------------------------------------------------------

We first exhibit proper submodules on the two singular hyperplanes.

\begin{prop}
\label{prop:reducibility}
Let \(P\) be a nonzero \(\mathcal D\)-module. If
\[
  (2\alpha-\beta)(2\alpha+\beta)=0,
\]
then \(F_{\alpha,\beta}(P)\) is reducible.
\end{prop}

\begin{proof}
Assume first that \(\mu_-=0\). We claim that \(P\oplus P\theta_1\) is a proper
submodule. It is proper because \(P\ne0\), so \(P\theta_2\ne0\). By
Lemma~\ref{lem:odd-component-formulas},
\[
  G_m^+P=0,\quad
  G_m^+(P\theta_1)\subseteq P,\quad
  G_m^-P\subseteq P\theta_1,\quad
  G_m^-(P\theta_1)=0,
\]
so \(P\oplus P\theta_1\) is stable under all odd generators. For the even
generators, the only terms in \(L_m\) and \(H_m\) that can send \(P\) into
\(P\Theta\) are
\[
  \frac{m^2}{2}\mu_-t^m\Theta
  \quad\text{and}\quad
  -m\mu_-t^m\Theta,
\]
both of which vanish when \(\mu_-=0\); the remaining even terms preserve
\(P\oplus P\theta_1\). Hence \(P\oplus P\theta_1\) is a nonzero proper
submodule.

Assume next that \(\mu_+=0\), and put
\[
  W=\operatorname{span}_{\mathbb C}
  \{\,v+2Dv\,\Theta,\ v\theta_2\mid v\in P\,\}.
\]
This subspace is nonzero and proper, since every vector in \(W\) has zero
\(P\theta_1\)-component while \(P\theta_1\ne0\).

We prove that \(W\) is stable. Let \(v\in P\). Since \(\mu_+=0\),
Lemma~\ref{lem:odd-component-formulas} gives
\[
  G_m^-v=-2t^mDv\,\theta_1,
  \quad
  G_m^-(Dv\,\Theta)=t^mDv\,\theta_1,
\]
so \(G_m^-(v+2Dv\,\Theta)=0\). Also, using \(D(t^m v)=t^m(D+m)v\),
\[
  G_m^-(v\theta_2)
  =-t^m v-2t^m(D+m)v\,\Theta
  =-\bigl(t^m v+2D(t^m v)\Theta\bigr)\in W .
\]
Hence \(W\) is stable under all \(G_m^-\). For \(G_m^+\) we have
\(G_m^+(v\theta_2)=0\) and
\[
  G_m^+(v+2Dv\,\Theta)
  =m\mu_-t^m v\,\theta_2+2t^mDv\,\theta_2
  =t^m(2D+m\mu_-)v\,\theta_2\in W ,
\]
so \(W\) is stable under all odd generators.

It remains to prove stability under the even generators. From the mixed odd
brackets
\[
  \{G_0^+,G_k^-\}=-2L_k-kH_k,
  \quad
  \{G_1^+,G_{k-1}^-\}=-2L_k+(2-k)H_k ,
\]
subtraction gives
\[
  H_k=
  \frac12\bigl(\{G_1^+,G_{k-1}^-\}-\{G_0^+,G_k^-\}\bigr),
  \quad
  L_k=-\frac12\{G_0^+,G_k^-\}-\frac{k}{2}H_k .
\]
Since \(W\) is stable under all odd generators, it is stable under the
right-hand sides, hence under all \(H_k\) and \(L_k\). Therefore \(W\) is a
nonzero proper submodule.
\end{proof}

The non-singular case is controlled by a second-difference operator that kills
the three higher exterior components \(P\theta_1\), \(P\theta_2\), \(P\Theta\)
and sends the \(P\)-component to \(P\Theta\).

\begin{lemm}
\label{lem:filter-operator}
Let \(k,d\in\mathbb Z\) with \(d\ne0\), and define
\[
  \mathcal S_{k,d}
  =
  \frac12\left(
    G_{-d}^+G_{k+d}^-
    +G_d^+G_{k-d}^-
    -2G_0^+G_k^-
  \right)
  \in U(\mathfrak g).
\]
Then \(\mathcal S_{k,d}\) annihilates \(P\theta_1\oplus P\theta_2\oplus
P\Theta\), and for
\[
  w=v_0+v_1\theta_1+v_2\theta_2+v_{12}\Theta
\]
one has
\[
  \mathcal S_{k,d}w
  =
  -d^2\mu_-\mu_+\,t^k v_0\,\Theta .
\]
\end{lemm}

\begin{proof}
We compute \(G_a^+G_b^-\) with \(a+b=k\) and then take the centered second
difference at \(a=0,\pm d\).

For \(v\in P\),
\[
  G_b^-v=-t^b(2D+b\mu_+)v\,\theta_1,
\]
and applying \(G_a^+\) gives
\[
  G_a^+G_b^-v
  =
  -t^k(2D+b\mu_+)v
  +a\mu_-\,t^k(2D+b\mu_+)v\,\Theta .
\]
Since \(b=k-a\), the scalar term is affine linear in \(a\), so its centered
second difference vanishes. For the \(\Theta\)-term, substituting
\((a,b)=(-d,k+d),(d,k-d),(0,k)\) gives
\[
\begin{aligned}
  \mathcal S_{k,d}v
  &=
  \frac12\mu_-
  \left(
    -d\bigl(2D+(k+d)\mu_+\bigr)
    +d\bigl(2D+(k-d)\mu_+\bigr)
  \right)t^k v\,\Theta       \\
  &=
  -d^2\mu_-\mu_+\,t^k v\,\Theta .
\end{aligned}
\]

For the remaining components, first \(G_b^-(P\theta_1)=0\), so
\(\mathcal S_{k,d}(P\theta_1)=0\). Next,
\[
  G_a^+G_b^-(v\Theta)=t^k v-a\mu_-t^k v\,\Theta
\]
is affine linear in \(a\), so its centered second difference vanishes and
\(\mathcal S_{k,d}(P\Theta)=0\). Finally,
\[
  G_b^-(v\theta_2)
  =
  -t^b v-t^b(2D+2b+b\mu_+)v\,\Theta ,
\]
and applying \(G_a^+\) gives
\[
\begin{aligned}
  G_a^+G_b^-(v\theta_2)
  &=
  -a\mu_-t^k v\,\theta_2
  -t^k(2D+b(\mu_++2))v\,\theta_2  \\
  &=
  -t^k\bigl(2D+k(\mu_++2)
      +a(\mu_- -\mu_+-2)\bigr)v\,\theta_2 ,
\end{aligned}
\]
which is again affine linear in \(a\). Hence its centered second difference
vanishes, and the lemma follows.
\end{proof}

\begin{theo}
\label{thm:irreducibility}
Let \(P\) be a simple \(\mathcal D\)-module. Then
\[
  F_{\alpha,\beta}(P)\ \text{is simple}
  \quad\Longleftrightarrow\quad
  (2\alpha-\beta)(2\alpha+\beta)\ne0 .
\]
\end{theo}

\begin{proof}
Necessity follows from Proposition~\ref{prop:reducibility}. Assume now that
\(\mu_-\mu_+\ne0\), and let \(0\ne W\subseteq F_{\alpha,\beta}(P)\) be a
\(\mathfrak g\)-submodule; we prove that \(W=F_{\alpha,\beta}(P)\).

Choose \(0\ne w=v_0+v_1\theta_1+v_2\theta_2+v_{12}\Theta\in W\). We first show
that \(W\) contains a nonzero vector in \(P\Theta\). If \(v_0\ne0\), then
Lemma~\ref{lem:filter-operator} gives
\[
  \mathcal S_{0,1}w=-\mu_-\mu_+v_0\Theta\in W,
\]
which is nonzero. If \(v_0=0\) but \(v_1\ne0\), then \(G_0^+w\) has nonzero
\(P\)-component. If \(v_0=v_1=0\) but \(v_2\ne0\), then \(G_0^-w\) has nonzero
\(P\)-component. In both cases we reduce to the previous one. If
\(v_0=v_1=v_2=0\), then \(w\in P\Theta\) already. Hence \(u\Theta\in W\) for
some \(0\ne u\in P\).

The zero-mode formulas
\[
  G_0^-(u\Theta)=u\theta_1,\quad
  G_0^+(u\Theta)=u\theta_2,\quad
  G_0^+(u\theta_1)=u
\]
imply \(u\otimes\Lambda(\theta_1,\theta_2)\subseteq W\). Define
\[
  P_W=
  \{\,p\in P\mid p\otimes\Lambda(\theta_1,\theta_2)\subseteq W\,\},
\]
so that \(P_W\ne0\). We prove that \(P_W\) is a \(\mathcal D\)-submodule of
\(P\).

Let \(p\in P_W\). Since \(p\Theta\in W\), for every \(m\in\mathbb Z\) we have
\[
  G_m^-(p\Theta)=t^m p\,\theta_1\in W,
  \quad
  G_m^+(p\Theta)=t^m p\,\theta_2\in W.
\]
Applying \(G_0^+\) to \(t^mp\,\theta_1\) gives \(t^mp\in W\), and applying
Lemma~\ref{lem:filter-operator} to this vector gives
\(-\mu_-\mu_+\,t^mp\,\Theta\in W\). Since \(\mu_-\mu_+\ne0\), we get
\(t^mp\,\Theta\in W\), and the zero-mode formulas then give
\(t^mp\otimes\Lambda(\theta_1,\theta_2)\subseteq W\). Thus \(t^mp\in P_W\) for
all \(m\in\mathbb Z\).

It remains to check stability under \(D\). Since \(L_0=D\) on the coefficient
factor, \(Dp=L_0p\in W\). If \(Dp\ne0\), applying
Lemma~\ref{lem:filter-operator} to \(Dp\) gives
\(-\mu_-\mu_+\,Dp\,\Theta\in W\), and the zero-mode formulas give
\(Dp\otimes\Lambda(\theta_1,\theta_2)\subseteq W\); if \(Dp=0\), this is
trivial. Hence \(Dp\in P_W\), and \(P_W\) is a \(\mathcal D\)-submodule of
\(P\).

Since \(P\) is simple and \(P_W\ne0\), we have \(P_W=P\), and therefore
\(W=F_{\alpha,\beta}(P)\). This proves simplicity in the non-singular case.
\end{proof}

%----------------------------------------------------------------------
\subsection{Generic isomorphism classes}
\label{subsec:isomorphism}
%----------------------------------------------------------------------

All isomorphisms in this subsection are even homomorphisms of
\(\mathfrak g\)-modules.

\begin{theo}[Generic isomorphism classes]
\label{thm:isomorphism}
Let \(P_1,P_2\) be simple \(\mathcal D\)-modules, and assume
\[
  (2\alpha_i-\beta_i)(2\alpha_i+\beta_i)\ne0
  \quad(i=1,2).
\]
Then \(F_{\alpha_1,\beta_1}(P_1)\cong F_{\alpha_2,\beta_2}(P_2)\) as
\(\mathfrak g\)-modules if and only if
\(\alpha_1=\alpha_2\), \(\beta_1=\beta_2\), and \(P_1\cong P_2\) as
\(\mathcal D\)-modules.
\end{theo}

\begin{proof}
Sufficiency is immediate: if \(\phi:P_1\to P_2\) is a \(\mathcal D\)-module
isomorphism and the parameters agree, then
\[
  v\xi\longmapsto \phi(v)\xi,
  \quad \xi\in\Lambda(\theta_1,\theta_2),
\]
intertwines the defining actions.

Conversely, let
\[
  \Phi:F_{\alpha_1,\beta_1}(P_1)
  \longrightarrow
  F_{\alpha_2,\beta_2}(P_2)
\]
be a \(\mathfrak g\)-module isomorphism, and put
\(\mu_{-,i}=2\alpha_i-\beta_i\), \(\mu_{+,i}=2\alpha_i+\beta_i\). On
\(F_{\alpha_i,\beta_i}(P_i)\), the operator \(H_0\) has the three distinct
eigenvalues \(\beta_i-1\), \(\beta_i\), \(\beta_i+1\) on \(P_i\theta_1\),
\(P_i\oplus P_i\Theta\), \(P_i\theta_2\), respectively.
Since \(\Phi\) commutes with \(H_0\), the unordered sets
\[
  \{\beta_1-1,\beta_1,\beta_1+1\},
  \quad
  \{\beta_2-1,\beta_2,\beta_2+1\}
\]
coincide; comparing their sums gives \(\beta_1=\beta_2\). Write the common
value as \(\beta\). Then \(\Phi\) preserves the three \(H_0\)-eigenspaces:
\[
  \Phi(P_1\theta_1)=P_2\theta_1,\quad
  \Phi(P_1\theta_2)=P_2\theta_2,\quad
  \Phi(P_1\oplus P_1\Theta)=P_2\oplus P_2\Theta.
\]

Inside the middle eigenspace, the operator \(\mathcal S_{0,1}\) separates
\(P_i\Theta\): by Lemma~\ref{lem:filter-operator},
\[
  \mathcal S_{0,1}(v_0+v_{12}\Theta)
  =
  -\mu_{-,i}\mu_{+,i}\,v_0\Theta,
\]
and since \(\mu_{-,i}\mu_{+,i}\ne0\), the kernel of \(\mathcal S_{0,1}\) on the
middle eigenspace is precisely \(P_i\Theta\). As \(\mathcal S_{0,1}\in
U(\mathfrak g)\), the map \(\Phi\) commutes with it, so
\(\Phi(P_1\Theta)=P_2\Theta\). Then
\[
  \Phi(P_1)=\Phi(G_0^+(P_1\theta_1))
  =G_0^+\Phi(P_1\theta_1)
  =G_0^+(P_2\theta_1)=P_2.
\]

Thus there exist linear isomorphisms
\(\phi,\phi_1,\phi_2,\phi_{12}:P_1\to P_2\) with
\[
  \Phi(v)=\phi(v),\quad
  \Phi(v\theta_1)=\phi_1(v)\theta_1,\quad
  \Phi(v\theta_2)=\phi_2(v)\theta_2,\quad
  \Phi(v\Theta)=\phi_{12}(v)\Theta .
\]
Commuting with
\[
  G_0^+(v\theta_1)=v,\quad
  G_0^+(v\Theta)=v\theta_2,\quad
  G_0^-(v\Theta)=v\theta_1
\]
gives \(\phi=\phi_1=\phi_2=\phi_{12}\); denote the common map by \(\phi\).

Now \(G_0^-(v\theta_2)=-v-2Dv\,\Theta\), so commuting \(\Phi\) with \(G_0^-\)
gives
\[
  -\phi(v)-2\phi(Dv)\Theta
  =
  -\phi(v)-2D\phi(v)\Theta ,
\]
hence \(\phi(Dv)=D\phi(v)\).

Next, for arbitrary \(m\in\mathbb Z\),
\[
  G_m^+(v\theta_1)=t^m v-m(2\alpha_1-\beta)t^m v\,\Theta ,
\]
so applying \(\Phi\) gives \(\phi(t^m v)-m(2\alpha_1-\beta)\phi(t^m v)\Theta\),
while
\[
  G_m^+(\Phi(v\theta_1))
  =
  t^m\phi(v)-m(2\alpha_2-\beta)t^m\phi(v)\Theta .
\]
Comparing \(P_2\)-components gives \(\phi(t^m v)=t^m\phi(v)\). Choose
\(m\ne0\) and \(0\ne v\in P_1\). Since \(t^m\) is invertible on every
\(\mathcal D\)-module and \(\phi\) is an isomorphism, we have
\(t^m\phi(v)\ne0\); comparing \(P_2\Theta\)-components then gives
\(2\alpha_1-\beta=2\alpha_2-\beta\), so \(\alpha_1=\alpha_2\). Thus \(\phi\)
commutes with \(D\) and with every \(t^m\), and hence
\(\phi:P_1\to P_2\) is a \(\mathcal D\)-module isomorphism.
\end{proof}

%----------------------------------------------------------------------
\subsection{The singular hyperplanes}
\label{subsec:singular-hyperplanes}
%----------------------------------------------------------------------

On the singular hyperplanes, the four-component modules become extensions of
two-component modules. The modules \(\mathcal T_b(P)\) below are the
two-component modules used in \cite{CDLP2}. We recall their formulas and also
record the charge-reversed version needed for one of the two singular
hyperplanes.

For a \(\mathcal D\)-module \(P\) and \(b\in\mathbb C\), let
\(\mathcal T_b(P)=P\oplus P\theta\) be defined by
\[
\begin{array}{ll}
  L_m v=t^m(D+bm)v,
&
  L_m(v\theta)=t^m\left(D+\left(b+\frac12\right)m\right)v\theta,
\\[3pt]
  H_m v=-2b\,t^m v,
&
  H_m(v\theta)=(1-2b)t^m v\theta,
\\[3pt]
  G_m^+v=2t^m(D+2bm)v\theta,
&
  G_m^+(v\theta)=0,
\\[3pt]
  G_m^-v=0,
&
  G_m^-(v\theta)=-t^m v .
\end{array}
\]
The charge-reversed two-component module
\(\mathcal T_b^{\mathrm r}(P)=P\oplus P\theta\) is defined by
\[
\begin{array}{ll}
  L_m v=t^m(D+bm)v,
&
  L_m(v\theta)=t^m\left(D+\left(b+\frac12\right)m\right)v\theta,
\\[3pt]
  H_m v=2b\,t^m v,
&
  H_m(v\theta)=(2b-1)t^m v\theta,
\\[3pt]
  G_m^+v=0,
&
  G_m^+(v\theta)=t^m v,
\\[3pt]
  G_m^-v=-2t^m(D+2bm)v\theta,
&
  G_m^-(v\theta)=0 .
\end{array}
\]
If \(N=N_{\bar0}\oplus N_{\bar1}\) is a \(\mathbb Z_2\)-graded module, then
\(\Pi N\) denotes the same linear \(\mathfrak g\)-action on the parity-reversed
superspace, with \((\Pi N)_{\bar0}=N_{\bar1}\) and \((\Pi N)_{\bar1}=N_{\bar0}\).
For ungraded modules the symbol \(\Pi\) may be omitted.

\begin{prop}
\label{prop:singular-exact-sequences}
Let \(P\) be a \(\mathcal D\)-module.  Then there are exact sequences of
\(\mathfrak g\)-modules
\[
  0\longrightarrow
  \mathcal T_{\alpha}^{\mathrm r}(P)
  \longrightarrow
  F_{\alpha,\,2\alpha}(P)
  \longrightarrow
  \Pi\mathcal T_{\alpha+\frac12}^{\mathrm r}(P)
  \longrightarrow0,
\]
and
\[
  0\longrightarrow
  \mathcal T_{\alpha}(P)
  \longrightarrow
  F_{\alpha,\,-2\alpha}(P)
  \longrightarrow
  \Pi\mathcal T_{\alpha+\frac12}(P)
  \longrightarrow0 .
\]
\end{prop}

\begin{proof}
Assume first that \(\beta=2\alpha\), so \(\mu_-=0\). The subspace
\(P\oplus P\theta_1\) is the stable submodule from
Proposition~\ref{prop:reducibility}, and on it
\[
  L_mv=t^m(D+\alpha m)v,
  \quad
  L_m(v\theta_1)=t^m\left(D+\left(\alpha+\frac12\right)m\right)v\theta_1,
\]
\[
  H_m v=2\alpha\,t^m v,
  \quad
  H_m(v\theta_1)=(2\alpha-1)t^m v\theta_1,
\]
\[
  G_m^+(v\theta_1)=t^m v,
  \quad
  G_m^-v=-2t^m(D+2\alpha m)v\,\theta_1,
  \quad
  G_m^+v=0,\quad G_m^-(v\theta_1)=0.
\]
Therefore \(P\oplus P\theta_1\cong\mathcal T_\alpha^{\mathrm r}(P)\).

The quotient is represented by \(P\theta_2\oplus P\Theta\); write
\(\overline{w}\) for the quotient class of \(w\). The even action is
\[
  L_m\overline{v\theta_2}
  =
  \overline{t^m\left(D+\left(\alpha+\frac12\right)m\right)v\,\theta_2},
  \quad
  L_m\overline{v\Theta}
  =
  \overline{t^m(D+(\alpha+1)m)v\,\Theta},
\]
\[
  H_m\overline{v\theta_2}
  =
  (2\alpha+1)\overline{t^m v\theta_2},
  \quad
  H_m\overline{v\Theta}
  =
  2\alpha \overline{t^mv\Theta},
\]
and the odd action is
\[
  G_m^+\overline{v\Theta}
  =
  \overline{t^m v\,\theta_2},
  \quad
  G_m^+\overline{v\theta_2}=0,
\]
\[
  G_m^-\overline{v\theta_2}
  =
  -2\overline{t^m(D+(2\alpha+1)m)v\,\Theta},
  \quad
  G_m^-\overline{v\Theta}=0.
\]
These are the formulas for \(\mathcal T_{\alpha+\frac12}^{\mathrm r}(P)\).
Since the quotient representatives \(P\theta_2\) and \(P\Theta\) have the
opposite parity from the standard \(P\oplus P\theta\) notation, the quotient
is \(\Pi\mathcal T_{\alpha+\frac12}^{\mathrm r}(P)\). This proves the first
exact sequence.

Assume now that \(\beta=-2\alpha\), so \(\mu_+=0\). Put
\[
  X(v)=v+2Dv\,\Theta,
  \quad
  Y(v)=v\theta_2 .
\]
The subspace
\(W=\operatorname{span}_{\mathbb C}\{\,X(v),Y(v)\mid v\in P\,\}\) is the stable
submodule from Proposition~\ref{prop:reducibility}; we identify \(X(v)\) with
the even component and \(Y(v)\) with the odd component. Direct calculation
gives
\[
  L_mX(v)=X\bigl(t^m(D+\alpha m)v\bigr),
  \quad
  L_mY(v)=Y\bigl(t^m(D+(\alpha+\tfrac12)m)v\bigr),
\]
\[
  H_mX(v)=-2\alpha\,X(t^m v),
  \quad
  H_mY(v)=(1-2\alpha)Y(t^m v),
\]
\[
  G_m^-X(v)=0,\quad
  G_m^-Y(v)=-X(t^m v),
  \quad
  G_m^+Y(v)=0,
  \quad
  G_m^+X(v)=2Y\bigl(t^m(D+2\alpha m)v\bigr),
\]
so \(W\cong\mathcal T_\alpha(P)\).

It remains to identify the quotient, represented by \(P\theta_1\oplus
P\Theta\). Define quotient vectors
\[
  U(v)=-\overline{v\theta_1},
  \quad
  V(v)=\overline{v\Theta}.
\]
The even action is
\[
  L_mU(v)=U\bigl(t^m(D+(\alpha+\tfrac12)m)v\bigr),
  \quad
  L_mV(v)=V\bigl(t^m(D+(\alpha+1)m)v\bigr),
\]
\[
  H_mU(v)=-(2\alpha+1)t^mU(v),
  \quad
  H_mV(v)=-2\alpha\,t^mV(v),
\]
which are the \(H\)-weights \(-2b\) and \(1-2b\) for \(b=\alpha+\frac12\). For
the odd action,
\[
  G_m^-V(v)=\overline{t^m v\,\theta_1}=-U(t^m v),
  \quad
  G_m^-U(v)=0,
\]
and \(G_m^+V(v)=\overline{t^m v\,\theta_2}=0\) in the quotient, since
\(t^m v\,\theta_2=Y(t^m v)\in W\). Finally,
\[
  G_m^+U(v)
  =
  -\overline{G_m^+(v\theta_1)}
  =
  -\overline{t^m v-4\alpha m\,t^m v\,\Theta}.
\]
Since \(X(t^m v)=t^m v+2D(t^m v)\Theta\in W\), in the quotient
\[
  \overline{t^m v}
  =
  -2\overline{D(t^m v)\Theta}
  =
  -2\overline{t^m(D+m)v\,\Theta},
\]
and therefore
\[
\begin{aligned}
  G_m^+U(v)
  &=
  2\overline{t^m(D+m)v\,\Theta}
  +4\alpha m\,\overline{t^m v\,\Theta}        \\
  &=
  2V\bigl(t^m(D+(2\alpha+1)m)v\bigr).
\end{aligned}
\]
These are the formulas for \(\mathcal T_{\alpha+\frac12}(P)\). Again the
natural parity of the quotient representatives is opposite to the standard
\(P\oplus P\theta\) notation, so the quotient is
\(\Pi\mathcal T_{\alpha+\frac12}(P)\). This proves the second exact sequence.
\end{proof}

%======================================================================
\section{Applications and examples}
\label{sec:applications}
%======================================================================

In this final section we apply the construction to two standard families of
simple \(\mathcal D\)-modules, obtaining non-weight and weight
\(\mathfrak g\)-modules, and we relate the two through Nilsson's weighting
functor.

We end with concrete choices of the underlying \(\mathcal D\)-module \(P\). The
construction applies to every simple \(\mathcal D\)-module; we spell out two
standard classes. The polynomial Weyl modules \(\Omega(\lambda)\),
\(\lambda\in\mathbb C^*\), give non-weight \(\mathfrak g\)-modules that are free
of rank \(4\) over \(U(\mathbb C L_0)\), while the rank-one weight Weyl modules
\(\mathcal V_a\) give bounded weight modules with four-dimensional weight
spaces. Nilsson's weighting functor with respect to \(L_0\) relates the two
families.

We keep the notation \(\mu_-=2\alpha-\beta\), \(\mu_+=2\alpha+\beta\), and,
unless otherwise stated, work in the non-singular region \(\mu_-\mu_+\ne0\). In
this region Theorem~\ref{thm:irreducibility} implies that
\(F_{\alpha,\beta}(P)\) is simple whenever \(P\) is a simple
\(\mathcal D\)-module.

%----------------------------------------------------------------------
\subsection{Rank-four free modules}
\label{subsec:free_modules}
%----------------------------------------------------------------------

Let \(\lambda\in\mathbb C^*\). Recall that \(\Omega(\lambda)=\mathbb C[x]\) is
the \(\mathcal D\)-module with
\[
  t^m f(x)=\lambda^m f(x-m),
  \quad
  Df(x)=xf(x),
  \quad m\in\mathbb Z .
\]

\begin{lemm}
\label{lem:Omega-irreducible}
For every \(\lambda\in\mathbb C^*\), the \(\mathcal D\)-module
\(\Omega(\lambda)\) is simple, and
\[
  \Omega(\lambda_1)\cong\Omega(\lambda_2)
  \quad\Longleftrightarrow\quad
  \lambda_1=\lambda_2 .
\]
\end{lemm}

\begin{proof}
Let \(0\ne W\subseteq\Omega(\lambda)\) be a \(\mathcal D\)-submodule, and choose
a nonzero \(f(x)\in W\) of minimal degree. Since \(W\) is stable under \(t\),
it is stable under \(\lambda^{-1}t-1\), so \(f(x-1)-f(x)\in W\). If
\(\deg f>0\), then \(f(x-1)-f(x)\) is nonzero of smaller degree than \(f\), a
contradiction. Thus \(W\) contains a nonzero constant. Since \(D\) acts as
multiplication by \(x\), it follows that \(W\) contains all powers of \(x\),
and hence \(W=\mathbb C[x]\).

If \(\phi:\Omega(\lambda_1)\to\Omega(\lambda_2)\) is a \(\mathcal D\)-module
isomorphism with \(\phi(1)=f(x)\ne0\), then the \(t\)-action gives
\(\lambda_2 f(x-1)=\lambda_1 f(x)\), and comparing leading coefficients gives
\(\lambda_1=\lambda_2\). The converse is immediate.
\end{proof}

Applying the construction of Section~\ref{sec:twisting-functors}, we set
\[
  \mathcal M_{\mathrm{free}}(\lambda;\alpha,\beta)
  :=
  F_{\alpha,\beta}(\Omega(\lambda)),
\]
so that, as a vector space,
\[
  \mathcal M_{\mathrm{free}}(\lambda;\alpha,\beta)
  =
  \mathbb C[x]\oplus\mathbb C[x]\theta_1
  \oplus\mathbb C[x]\theta_2
  \oplus\mathbb C[x]\Theta .
\]
Since \(L_0=D\) acts as multiplication by \(x\), the elements
\(1,\theta_1,\theta_2,\Theta\) form a free basis over
\(U(\mathbb C L_0)\cong\mathbb C[L_0]\). Thus
\(\mathcal M_{\mathrm{free}}(\lambda;\alpha,\beta)\) is \(U(\mathbb C L_0)\)-free
of rank \(4\).

We write elements as column vectors with respect to the ordered exterior basis
\(1,\theta_1,\theta_2,\Theta\):
\[
  f_0(x)+f_1(x)\theta_1+f_2(x)\theta_2+f_{12}(x)\Theta
  \longleftrightarrow
  \mathbf f(x)=
  \begin{pmatrix}
    f_0(x)\\
    f_1(x)\\
    f_2(x)\\
    f_{12}(x)
  \end{pmatrix},
\]
and let \(S_{-m}\) be the shift operator \((S_{-m}\mathbf f)(x)=\mathbf
f(x-m)\).

\begin{prop}
\label{prop:free-module-formulas}
For \(m\in\mathbb Z\), the action of \(\mathfrak g\) on
\(\mathcal M_{\mathrm{free}}(\lambda;\alpha,\beta)\cong\mathbb C[x]^4\) is
given by
\[
  L_m\mathbf f(x)=\lambda^m M_L(m,x)S_{-m}\mathbf f(x),
  \quad
  H_m\mathbf f(x)=\lambda^m M_H(m,x)S_{-m}\mathbf f(x),
\]
\[
  G_m^\pm\mathbf f(x)=\lambda^m M_{G^\pm}(m,x)S_{-m}\mathbf f(x),
\]
where
\[
  M_L(m,x)=
  \begin{pmatrix}
    x-m+m\alpha & 0 & 0 & 0 \\
    0 & x-m+m(\alpha+\frac12) & 0 & 0 \\
    0 & 0 & x-m+m(\alpha+\frac12) & 0 \\
    \frac{m^2}{2}\mu_- & 0 & 0 & x-m+m(\alpha+1)
  \end{pmatrix},
\]
\[
  M_H(m,x)=
  \begin{pmatrix}
    \beta & 0 & 0 & 0 \\
    0 & \beta-1 & 0 & 0 \\
    0 & 0 & \beta+1 & 0 \\
    -m\mu_- & 0 & 0 & \beta
  \end{pmatrix},
\]
\[
  M_{G^+}(m,x)=
  \begin{pmatrix}
    0 & 1 & 0 & 0 \\
    0 & 0 & 0 & 0 \\
    m\mu_- & 0 & 0 & 1 \\
    0 & -m\mu_- & 0 & 0
  \end{pmatrix},
\]
and
\[
  M_{G^-}(m,x)=
  \begin{pmatrix}
    0 & 0 & -1 & 0 \\
    -2\bigl(x-m+\frac m2\mu_+\bigr) & 0 & 0 & 1 \\
    0 & 0 & 0 & 0 \\
    0 & 0 & -2\bigl(x+\frac m2\mu_+\bigr) & 0
  \end{pmatrix}.
\]
\end{prop}

\begin{proof}
Substitute the \(\Omega(\lambda)\)-action into the component formulas for
\(F_{\alpha,\beta}(P)\); the only rules needed are
\(t^m f(x)=\lambda^m f(x-m)\) and \(t^mD f(x)=\lambda^m(x-m)f(x-m)\). For
example,
\[
  t^m(D+m\alpha)f(x)
  =
  \lambda^m\bigl(x-m+m\alpha\bigr)f(x-m),
\]
which gives the \((1,1)\)-entry of \(M_L(m,x)\), while
\[
  -t^m(2D+m\mu_+)f(x)
  =
  -2\lambda^m\bigl(x-m+\tfrac m2\mu_+\bigr)f(x-m),
\]
which gives the \((2,1)\)-entry of \(M_{G^-}(m,x)\). The remaining entries
follow in the same way from Lemmas~\ref{lem:odd-component-formulas}
and~\ref{lem:even-component-formulas}.
\end{proof}

\begin{coro}
\label{cor:free-module-classification}
The module \(\mathcal M_{\mathrm{free}}(\lambda;\alpha,\beta)\) is simple if
and only if \(\mu_-\mu_+\ne0\). Moreover, in the non-singular region,
\[
  \mathcal M_{\mathrm{free}}(\lambda_1;\alpha_1,\beta_1)
  \cong
  \mathcal M_{\mathrm{free}}(\lambda_2;\alpha_2,\beta_2)
\]
as \(\mathfrak g\)-modules if and only if
\[
  \lambda_1=\lambda_2,\quad
  \alpha_1=\alpha_2,\quad
  \beta_1=\beta_2 .
\]
\end{coro}

\begin{proof}
This follows from Theorems~\ref{thm:irreducibility}
and~\ref{thm:isomorphism}, together with Lemma~\ref{lem:Omega-irreducible}.
\end{proof}

\begin{rema}
The module \(\mathcal M_{\mathrm{free}}(\lambda;\alpha,\beta)\) is not a weight
module for \(L_0\): indeed, \(L_0\) acts as multiplication by \(x\) on
\(\mathbb C[x]\), which has no nonzero eigenvectors. Thus these examples are
non-weight modules, even though they are free of finite rank over
\(U(\mathbb C L_0)\).
\end{rema}

\paragraph{A characterization of the free modules.}
\label{sec:reconstruction-free-modules}

We now give a converse description of the polynomial free modules
\(\mathcal M_{\mathrm{free}}(\lambda;\alpha,\beta)\): a recognition criterion
for the modules constructed above, rather than a classification of all
\(U(\mathbb C L_0)\)-free rank-four \(\mathfrak g\)-modules.

Put \(R=U(\mathbb C L_0)\cong\mathbb C[x]\), where \(x=L_0\), and let
\[
  \mathfrak g^{(0)}
  =
  \operatorname{span}_{\mathbb C}\{L_0,H_0,G_0^+,G_0^-\}
\]
be the degree-zero Lie subsuperalgebra of \(\mathfrak g\). Under the
normalization
\[
  C=-2L_0,\quad H=H_0,\quad E=G_0^+,\quad F=G_0^-,
\]
its defining relations become
\[
  [C,\mathfrak g^{(0)}]=0,\quad
  [H,E]=E,\quad [H,F]=-F,\quad [E,F]=C,
  \quad [E,E]=[F,F]=0 ,
\]
so that \(\mathfrak g^{(0)}\cong\mathfrak{gl}(1|1)\).

We use an \(R\)-free version of the standard zero-layer module. It is not a
weight module for the full Cartan subalgebra
\(\mathbb C L_0\oplus\mathbb C H_0\), since \(L_0\) acts as multiplication by
\(x\); it is an \(H_0\)-weight module whose \(H_0\)-weight spaces are free
\(R\)-modules.

\begin{defi}
\label{def:standard-zero-layer-module}
For \(\beta\in\mathbb C\), let \(S_\beta\) be the \(R\)-free
\(\mathfrak g^{(0)}\)-module with \(R\)-basis \(v_{12},v_1,v_2,v_0\). The
element \(L_0\) acts as multiplication by \(x\), and
\[
  H_0v_{12}=\beta v_{12},\quad
  H_0v_0=\beta v_0,
  \quad
  H_0v_1=(\beta-1)v_1,\quad
  H_0v_2=(\beta+1)v_2.
\]
The odd zero modes act by
\[
  G_0^-v_{12}=v_1,\quad
  G_0^+v_{12}=v_2,\quad
  G_0^+v_1=v_0,\quad
  G_0^-v_2=-v_0-2xv_{12},
\]
and
\[
  G_0^+v_0=G_0^+v_2=0,\quad
  G_0^-v_1=0,\quad
  G_0^-v_0=-2xv_1 .
\]
The parity of \(v_{12}\) is arbitrary; changing it gives the parity shift.
\end{defi}

The formulas in Definition~\ref{def:standard-zero-layer-module} satisfy the
relations of \(\mathfrak g^{(0)}\). Indeed, \(v_1=G_0^-v_{12}\),
\(v_2=G_0^+v_{12}\), and \(v_0=G_0^+G_0^-v_{12}\), and the remaining actions
follow from
\[
  (G_0^+)^2=(G_0^-)^2=0,
  \quad
  G_0^+G_0^-+G_0^-G_0^+=-2L_0 .
\]

Let
\[
  \mathfrak h=
  \operatorname{span}_{\mathbb C}\{L_m,H_m\mid m\in\mathbb Z\}
\]
be the Heisenberg--Virasoro subalgebra of \(\mathfrak g\). We use the
following classification of rank-one \(R\)-free \(\mathfrak h\)-modules.

\begin{theo}[cf.~\cite{CG}]
\label{thm:rank-one-even-free}
Let \(N\) be an \(\mathfrak h\)-module that is free of rank \(1\) over
\(R=U(\mathbb C L_0)\). Then there exist \(\lambda\in\mathbb C^*\) and
\(a,b\in\mathbb C\) such that \(N\cong\Omega_{\bar0}(\lambda,a,b)\), where
\(\Omega_{\bar0}(\lambda,a,b)=\mathbb C[x]\omega\) and
\[
  L_m f(x)\omega
  =
  \lambda^m(x-ma)f(x-m)\omega,
  \quad
  H_m f(x)\omega
  =
  b\lambda^m f(x-m)\omega .
\]
Moreover, \(\Omega_{\bar0}(\lambda,a,b)\) is simple if and only if
\((a,b)\ne(0,0)\).
\end{theo}

\begin{defi}
\label{def:standard-rank-four-free}
A \(\mathfrak g\)-module \(M\) is a \emph{standard rank-four free module} if
there exist \(\beta\in\mathbb C\) and a homogeneous vector \(\omega\in M\)
such that
\begin{enumerate}
\item
\(M|_{\mathfrak g^{(0)}}\cong S_\beta\), with \(\omega\) corresponding to
\(v_{12}\); and

\item
\(R\omega\) is a rank-one \(R\)-free \(\mathfrak h\)-submodule of \(M\).
\end{enumerate}
\end{defi}

Thus, if \(M\) is standard rank-four free, then
\[
  M=
  R\omega\oplus RG_0^-\omega\oplus RG_0^+\omega
  \oplus RG_0^+G_0^-\omega
\]
as an \(R\)-module.

\begin{prop}
\label{prop:F-standard-rank-four}
For every \(\lambda\in\mathbb C^*\) and \(\alpha,\beta\in\mathbb C\), the
module \(F_{\alpha,\beta}(\Omega(\lambda))\) is a standard rank-four free
module; one may take \(\omega=\Theta\).
\end{prop}

\begin{proof}
Let \(M=F_{\alpha,\beta}(\Omega(\lambda))\). As an \(R=\mathbb C[L_0]\)-module,
\(M\) has the free basis \(\Theta,\theta_1,\theta_2,1\). The zero-mode formulas
give
\[
  G_0^-\Theta=\theta_1,\quad
  G_0^+\Theta=\theta_2,\quad
  G_0^+G_0^-\Theta=1,
\]
together with \(H_0\Theta=\beta\Theta\), \(H_0 1=\beta1\),
\(H_0\theta_1=(\beta-1)\theta_1\), and \(H_0\theta_2=(\beta+1)\theta_2\).
Hence \(M|_{\mathfrak g^{(0)}}\cong S_\beta\), with \(\Theta\) corresponding to
\(v_{12}\).

The even-part formulas are
\[
  L_m\Theta=\lambda^m(x+m\alpha)\Theta,
  \quad
  H_m\Theta=\beta\lambda^m\Theta .
\]
Thus \(R\Theta\) is a rank-one \(R\)-free \(\mathfrak h\)-submodule, with
\(R\Theta\cong\Omega_{\bar0}(\lambda,-\alpha,\beta)\), and \(M\) is standard
rank-four free.
\end{proof}

\begin{theo}
\label{thm:characterization-standard-rank-four}
Let \(M\) be a \(\mathfrak g\)-module. Then \(M\) is a standard rank-four free
module if and only if there exist \(\lambda\in\mathbb C^*\) and
\(\alpha,\beta\in\mathbb C\) with
\[
  M\cong F_{\alpha,\beta}(\Omega(\lambda))
\]
as \(\mathfrak g\)-modules, up to parity. More precisely, if \(M\) is standard
with generator \(\omega\) and
\[
  L_m\omega=\lambda^m(x-ma)\omega,
  \quad
  H_m\omega=\beta\lambda^m\omega ,
\]
then the corresponding parameter is \(\alpha=-a\).
\end{theo}

\begin{proof}
One direction is Proposition~\ref{prop:F-standard-rank-four}. Conversely,
assume that \(M\) is standard and choose \(\omega\) as in
Definition~\ref{def:standard-rank-four-free}. Since \(R\omega\) is a rank-one
\(R\)-free \(\mathfrak h\)-module, Theorem~\ref{thm:rank-one-even-free} gives
\[
  L_m\omega=\lambda^m(x-ma)\omega,
  \quad
  H_m\omega=\beta\lambda^m\omega
\]
for some \(\lambda\in\mathbb C^*\) and \(a,\beta\in\mathbb C\). Setting
\(\alpha=-a\), we have
\[
  L_m\omega=\lambda^m(x+m\alpha)\omega,
  \quad
  H_m\omega=\beta\lambda^m\omega .
\]

Put
\[
  e_{12}=\omega,\quad
  e_1=G_0^-\omega,\quad
  e_2=G_0^+\omega,\quad
  e_0=G_0^+G_0^-\omega .
\]
By the standard zero-layer condition,
\[
  M=Re_0\oplus Re_1\oplus Re_2\oplus Re_{12},
\]
and the \(H_0\)-eigenspaces are
\[
  M_{\beta-1}=Re_1,\quad
  M_{\beta+1}=Re_2,\quad
  M_\beta=Re_{12}\oplus Re_0 .
\]

Since \(L_m\) and \(H_m\) preserve \(H_0\)-weights, the \(R\)-lines \(Re_1\)
and \(Re_2\) are stable under \(\mathfrak h\). Applying
Theorem~\ref{thm:rank-one-even-free} to these two lines, write
\[
  L_m e_2=\lambda_+^m(x-ma_+)e_2,
  \quad
  H_m e_2=(\beta+1)\lambda_+^m e_2,
\]
and
\[
  L_m e_1=\lambda_-^m(x-ma_-)e_1,
  \quad
  H_m e_1=(\beta-1)\lambda_-^m e_1 .
\]

Write \(G_m^+e_{12}=p_m e_2\) and \(G_m^-e_{12}=q_m e_1\), with
\(p_m,q_m\in R\). Applying \([H_m,G_0^+]=G_m^+\) to \(e_{12}\) gives
\[
  p_m=(\beta+1)\lambda_+^m-\beta\lambda^m ,
\]
and applying \([L_m,G_0^+]=-\frac m2G_m^+\) to \(e_{12}\) gives
\[
  \lambda_+^m(x-ma_+)-\lambda^m(x+m\alpha)
  =
  -\frac m2p_m .
\]
Comparing the coefficient of \(x\) gives \(\lambda_+^m=\lambda^m\) for all
\(m\), hence \(\lambda_+=\lambda\). Thus \(p_m=\lambda^m\), and comparing
constant terms gives \(a_+=\frac12-\alpha\). Therefore
\[
  L_m e_2
  =
  \lambda^m\left(x+m\left(\alpha-\frac12\right)\right)e_2,
  \quad
  H_m e_2=(\beta+1)\lambda^m e_2,
  \quad
  G_m^+e_{12}=\lambda^m e_2 .
\]

The same argument, using \([H_m,G_0^-]=-G_m^-\) and
\([L_m,G_0^-]=-\frac m2G_m^-\), gives
\[
  L_m e_1
  =
  \lambda^m\left(x+m\left(\alpha-\frac12\right)\right)e_1,
  \quad
  H_m e_1=(\beta-1)\lambda^m e_1,
  \quad
  G_m^-e_{12}=\lambda^m e_1 .
\]

The zero-mode relation gives \(G_0^-e_2=-e_0-2xe_{12}\), and using
\(\{G_m^+,G_0^-\}=-2L_m+mH_m\) on \(e_{12}\) gives
\[
  G_m^+e_1
  =
  \lambda^m e_0
  -m(2\alpha-\beta)\lambda^m e_{12}.
\]
Applying \([L_m,G_0^+]=-\frac m2G_m^+\) and \([H_m,G_0^+]=G_m^+\) to \(e_1\),
we obtain
\[
  L_m e_0
  =
  \lambda^m(x+m(\alpha-1))e_0
  +
  \frac{m^2}{2}(2\alpha-\beta)\lambda^m e_{12},
\]
and
\[
  H_m e_0
  =
  \beta\lambda^m e_0
  -m(2\alpha-\beta)\lambda^m e_{12}.
\]
The relation \(\{G_m^+,G_0^+\}=0\) gives
\[
  G_m^+e_0=m(2\alpha-\beta)\lambda^m e_2,
  \quad
  G_m^+e_2=0 .
\]
Since \(G_m^-\) lowers the \(H_0\)-weight and there is no weight space
\(M_{\beta-2}\), one has \(G_m^-e_1=0\). Finally, using
\(\{G_0^+,G_m^-\}=-2L_m-mH_m\) first on \(e_1\) and then on \(e_{12}\), we
obtain
\[
  G_m^-e_0
  =
  -2\lambda^m
  \left(x-m+\frac m2(2\alpha+\beta)\right)e_1,
\]
and
\[
  G_m^-e_2
  =
  -\lambda^m e_0
  -2\lambda^m
  \left(x+\frac m2(2\alpha+\beta)\right)e_{12}.
\]

These are exactly the defining formulas of
\(F_{\alpha,\beta}(\Omega(\lambda))\) under the identification
\[
  e_0\mapsto 1,\quad
  e_1\mapsto \theta_1,\quad
  e_2\mapsto \theta_2,\quad
  e_{12}\mapsto \Theta .
\]
Since every homogeneous generator \(X_m\) satisfies \(X_m f(x)=f(x-m)X_m\) for
\(f(x)\in R\), the identification extends from the four generators to all of
\(M\). Hence \(M\cong F_{\alpha,\beta}(\Omega(\lambda))\). If \(\omega\) has the
opposite parity from \(\Theta\), the same argument gives the parity shift.
\end{proof}

\begin{rema}
Definition~\ref{def:standard-rank-four-free} is not meant to classify all
\(R\)-free rank-four \(\mathfrak g\)-modules; it is a recognition criterion for
the modules constructed in this paper. The degree-zero subalgebra
\(\mathfrak g^{(0)}\cong\mathfrak{gl}(1|1)\) fixes the four \(R\)-summands, and
the Heisenberg--Virasoro subalgebra \(\mathfrak h\) fixes the initial line
\(R\omega\).
\end{rema}

%----------------------------------------------------------------------
\subsection{Weight modules from rank-one \texorpdfstring{\(\mathcal D\)}
{D}-modules}
\label{subsec:rank-one-weight-modules}
%----------------------------------------------------------------------

For \(a\in\mathbb C\), let \(\mathcal V_a=\mathbb C[t^{\pm1}]\) be the rank-one
weight \(\mathcal D\)-module defined by
\[
  t^m\cdot t^k=t^{k+m},
  \quad
  D\cdot t^k=(k+a)t^k,
  \quad m,k\in\mathbb Z.
\]

\begin{lemm}
\label{lem:Va-irreducible}
For every \(a\in\mathbb C\), the \(\mathcal D\)-module \(\mathcal V_a\) is
simple, and
\[
  \mathcal V_a\cong\mathcal V_b
  \quad\Longleftrightarrow\quad
  a-b\in\mathbb Z.
\]
\end{lemm}

\begin{proof}
Let \(0\ne W\subseteq\mathcal V_a\) be a \(\mathcal D\)-submodule, and choose
\(0\ne w=\sum_{i=p}^q c_i t^i\in W\) with \(c_pc_q\ne0\). The operator \(D\)
acts diagonally with pairwise distinct eigenvalues \(a+i\), \(i\in\mathbb Z\),
so a polynomial in \(D\) can isolate a nonzero coefficient of \(w\); thus
\(t^j\in W\) for some \(j\). Applying \(t^m\) for \(m\in\mathbb Z\) gives
\(t^{j+m}\in W\) for all \(m\), and hence \(W=\mathcal V_a\).

If \(\mathcal V_a\cong\mathcal V_b\), their \(D\)-spectra coincide, so
\(a+\mathbb Z=b+\mathbb Z\), that is, \(a-b\in\mathbb Z\). Conversely, if
\(a-b=n\in\mathbb Z\), then \(t^k\mapsto t^{k+n}\) defines a
\(\mathcal D\)-module isomorphism \(\mathcal V_a\to\mathcal V_b\).
\end{proof}

Applying the construction gives the weight module
\(\mathcal I(a;\alpha,\beta):=F_{\alpha,\beta}(\mathcal V_a)\), or equivalently
\[
  \mathcal I(a;\alpha,\beta)
  =
  \bigoplus_{k\in\mathbb Z}
  \left(
  \mathbb C u_k\oplus\mathbb C p_k
  \oplus\mathbb C q_k\oplus\mathbb C r_k
  \right),
\]
where \(u_k=t^k\), \(p_k=t^k\theta_1\), \(q_k=t^k\theta_2\), and
\(r_k=t^k\Theta\). This is a bounded weight module with \(L_0\)-support
\(a+\mathbb Z\): each \(L_0\)-weight space is four-dimensional, and the four
exterior components carry the \(H_0\)-charges \(\beta,\beta-1,\beta+1,\beta\),
respectively.

\begin{prop}%[Action on the rank-one weight example]
\label{prop:intermediate-formulas}
Set \(X_k=(u_k,p_k,q_k,r_k)^{\mathsf T}\) and \(\ell_k=k+a\). Then
\[
  L_mX_k=
  \begin{pmatrix}
  \ell_k+m\alpha & 0 & 0 & 0\\
  0 & \ell_k+m(\alpha+\frac12) & 0 & 0\\
  0 & 0 & \ell_k+m(\alpha+\frac12) & 0\\
  \frac{m^2}{2}\mu_- & 0 & 0 & \ell_k+m(\alpha+1)
  \end{pmatrix}X_{k+m},
\]
\[
  H_mX_k=
  \begin{pmatrix}
  \beta & 0 & 0 & 0\\
  0 & \beta-1 & 0 & 0\\
  0 & 0 & \beta+1 & 0\\
  -m\mu_- & 0 & 0 & \beta
  \end{pmatrix}X_{k+m},
\]
\[
  G_m^+X_k=
  \begin{pmatrix}
  0 & 1 & 0 & 0\\
  0 & 0 & 0 & 0\\
  m\mu_- & 0 & 0 & 1\\
  0 & -m\mu_- & 0 & 0
  \end{pmatrix}X_{k+m},
\]
and
\[
  G_m^-X_k=
  \begin{pmatrix}
  0 & 0 & -1 & 0\\
  -2\bigl(\ell_k+\frac m2\mu_+\bigr) & 0 & 0 & 1\\
  0 & 0 & 0 & 0\\
  0 & 0 & -2\bigl(\ell_k+m+\frac m2\mu_+\bigr) & 0
  \end{pmatrix}X_{k+m}.
\]
The columns of these matrices give the actions on \(u_k,p_k,q_k,r_k\),
respectively.
\end{prop}

\begin{proof}
Substitute \(t^m t^k=t^{k+m}\) and \(Dt^k=(k+a)t^k=\ell_k t^k\) into the
component formulas for \(F_{\alpha,\beta}(P)\). For instance,
\[
  G_m^-u_k
  =
  -t^m(2D+m\mu_+)t^k\theta_1
  =
  -2\bigl(\ell_k+\tfrac m2\mu_+\bigr)p_{k+m},
\]
and
\[
  G_m^-q_k
  =
  -u_{k+m}
  -2\bigl(\ell_k+m+\tfrac m2\mu_+\bigr)r_{k+m}.
\]
The remaining entries follow from Lemmas~\ref{lem:odd-component-formulas}
and~\ref{lem:even-component-formulas}.
\end{proof}

\begin{coro}
\label{cor:intermediate-irreducibility}
The module \(\mathcal I(a;\alpha,\beta)\) is simple if and only if
\((2\alpha-\beta)(2\alpha+\beta)\ne0\). Moreover, if both parameter pairs are
non-singular, then
\(\mathcal I(a_1;\alpha_1,\beta_1)\cong\mathcal I(a_2;\alpha_2,\beta_2)\)
if and only if \(\alpha_1=\alpha_2\), \(\beta_1=\beta_2\), and
\(a_1-a_2\in\mathbb Z\).
\end{coro}

\begin{proof}
This follows from Theorems~\ref{thm:irreducibility}
and~\ref{thm:isomorphism}, together with Lemma~\ref{lem:Va-irreducible}.
\end{proof}

%----------------------------------------------------------------------
\subsection{Weighting the free modules}
\label{subsec:weighting-functor}
%----------------------------------------------------------------------

We conclude by relating the free and weight examples through Nilsson's
weighting functor~\cite{N2}; we need only its elementary form with respect to
the \(L_0\)-spectrum.

For a \(\mathcal D\)-module \(P\), set
\[
  \mathscr W_a^{\mathcal D}(P)
  :=
  \bigoplus_{k\in\mathbb Z} P/(D-a-k)P .
\]
The operator \(t^m\) sends the summand indexed by \(a+k\) to the one indexed by
\(a+k+m\). For the free module \(\Omega(\lambda)\), the map
\[
  [f(x)]_{a+k}
  \longmapsto
  f(a+k)\lambda^{-k}t^k
\]
gives a \(\mathcal D\)-module isomorphism
\[
  \mathscr W_a^{\mathcal D}(\Omega(\lambda))
  \cong
  \mathcal V_a .
\]

We apply the weighting functor to \(\mathfrak g\)-modules on which \(L_0\) acts
and the generators of \(\mathfrak g\) act homogeneously with respect to
\(\operatorname{ad}L_0\). For such a \(\mathfrak g\)-module \(M\) and
\(\delta\in\mathbb C\), define
\[
  \mathscr W_\delta(M)
  :=
  \bigoplus_{k\in\mathbb Z} M/\mathfrak m_{\delta+k}M,
  \quad
  \mathfrak m_\nu=(L_0-\nu)\subset U(\mathbb C L_0).
\]
If \(X_m\in\mathfrak g\) satisfies \([L_0,X_m]=mX_m\), then
\(X_m[v]_\nu=[X_mv]_{\nu+m}\); this is well-defined because
\[
  (L_0-(\nu+m))X_m=X_m(L_0-\nu).
\]

For the modules considered here, the two weighting constructions fit into the
commutative diagram
\[
\xymatrix@C=4.3em@R=3.2em{
  \Omega(\lambda)
    \ar[r]^-{F_{\alpha,\beta}}
    \ar[d]_-{\mathscr W_a^{\mathcal D}}
  &
  \mathcal M_{\mathrm{free}}(\lambda;\alpha,\beta)
    \ar[d]^-{\mathscr W_a}
  \\
  \mathcal V_a
    \ar[r]_-{F_{\alpha,\beta}}
  &
  \mathcal I(a;\alpha,\beta).
}
\]
The next theorem gives the explicit isomorphism on the right vertical side.

\begin{theo}
\label{thm:weighting-reduction}
For every \(a\in\mathbb C\), there is a natural isomorphism of
\(\mathfrak g\)-modules
\[
  \mathscr W_a(\mathcal M_{\mathrm{free}}(\lambda;\alpha,\beta))
  \cong
  \mathcal I(a;\alpha,\beta).
\]
In particular, the parameter \(\lambda\in\mathbb C^*\) disappears after
weighting.
\end{theo}

\begin{proof}
Set \(M=\mathcal M_{\mathrm{free}}(\lambda;\alpha,\beta)\) and
\(N=\mathcal I(a;\alpha,\beta)\), and write \(\nu_k=a+k\). Since \(L_0=D\) on
\(M\), quotienting by \(L_0-\nu_k\) amounts to imposing \(x=\nu_k\), so
\[
  M/\mathfrak m_{\nu_k}M
  \cong
  \bigl(\mathbb C[x]/(x-\nu_k)\bigr)
  \otimes\Lambda(\theta_1,\theta_2).
\]
Define
\[
  \Psi_k:
  M/\mathfrak m_{\nu_k}M
  \longrightarrow
  \mathbb C t^k\otimes\Lambda(\theta_1,\theta_2)
  \subset N
\]
by
\[
  \Psi_k\bigl([f(x)\xi]_{\nu_k}\bigr)
  =
  f(\nu_k)\lambda^{-k}t^k\xi,
  \quad
  f(x)\in\mathbb C[x],\quad
  \xi\in\Lambda(\theta_1,\theta_2).
\]
This is well-defined, since the kernel of evaluation at \(x=\nu_k\) is
\((x-\nu_k)\), and it is a linear isomorphism on each summand. Hence
\(\Psi=\bigoplus_{k\in\mathbb Z}\Psi_k\) is a vector-space isomorphism
\[
  \mathscr W_a(M)\longrightarrow N.
\]

It remains to prove that \(\Psi\) is \(\mathfrak g\)-equivariant. Let
\(\Xi=\{1,\theta_1,\theta_2,\Theta\}\). For
\(X_m\in\{L_m,H_m,G_m^+,G_m^-\}\) and \(\xi\in\Xi\),
Proposition~\ref{prop:free-module-formulas} gives
\[
  X_m(f(x)\xi)
  =
  \lambda^m
  \sum_{\eta\in\Xi}
  p^X_{\xi,\eta}(x,m)f(x-m)\eta,
\]
where the \(p^X_{\xi,\eta}(x,m)\) are polynomials in \(x\). In
\(N=\mathcal I(a;\alpha,\beta)\), the same coefficient polynomials are
evaluated at the target \(L_0\)-weight.

Take \([f(x)\xi]_{\nu_{k-m}}\) in the summand indexed by \(\nu_{k-m}\).
Applying \(\Psi_{k-m}\) first gives
\[
  \Psi_{k-m}\bigl([f(x)\xi]_{\nu_{k-m}}\bigr)
  =
  f(\nu_{k-m})\lambda^{-(k-m)}t^{k-m}\xi ,
\]
and then applying \(X_m\) in \(N\) gives
\[
  \sum_{\eta\in\Xi}
  p^X_{\xi,\eta}(\nu_k,m)
  f(\nu_k-m)\lambda^{-(k-m)}t^k\eta ,
\]
where \(\nu_k-m=\nu_{k-m}\). On the other hand, applying \(X_m\) first in the
weighted module and then \(\Psi_k\) gives
\[
\begin{aligned}
  \Psi_k\bigl([X_m(f(x)\xi)]_{\nu_k}\bigr)
  &=
  \Psi_k\left(
    \left[
    \lambda^m
    \sum_{\eta\in\Xi}
    p^X_{\xi,\eta}(x,m)f(x-m)\eta
    \right]_{\nu_k}
  \right)                                      \\
  &=
  \sum_{\eta\in\Xi}
  p^X_{\xi,\eta}(\nu_k,m)
  f(\nu_k-m)\lambda^{m-k}t^k\eta .
\end{aligned}
\]
Since \(\lambda^{m-k}=\lambda^{-(k-m)}\), the two expressions agree. Thus
\(\Psi\) intertwines all generators of \(\mathfrak g\) and is therefore a
\(\mathfrak g\)-module isomorphism.
\end{proof}

\end{document}